\input amstex\documentstyle{amsppt}  
\pagewidth{12.5cm}\pageheight{19cm}\magnification\magstep1
\topmatter
\title Hecke algebras and involutions in Weyl groups\endtitle
\author George Lusztig and David A. Vogan, Jr.\endauthor
\address{Department of Mathematics, M.I.T., Cambridge, MA 02139}\endaddress
\thanks{Supported in part by National Science Foundation grants DMS-0758262, DMS-0967272.}\endthanks
\endtopmatter   
\document
\define\Irr{\text{\rm Irr}}

\define\tph{\ti{\ph}}
\define\uca{\un{\ca}}
\define\ufH{\un{\fH}}

\define\bco{\bar{\co}}

\define\uM{\un M}

\define\dsv{\dashv}

\define\si{\sim}

\define\sqc{\sqcup}

\define\qua{\quad}

\define\bK{\bar K}

\define\op{\oplus}
   
\redefine\sp{\spadesuit}
\define\part{\partial}
\define\em{\emptyset}

\define\n{\notin}

\define\m{\mapsto}
\define\do{\dots}

\define\lra{\leftrightarrow}

\define\sub{\subset}    

\define\T{\times}
\define\ti{\tilde}
\define\nl{\newline}
\redefine\i{^{-1}}

\define\un{\underline}
\define\ov{\overline}
\define\ot{\otimes}
\define\bbq{\bar{\QQ}_l}
\define\ql{\QQ_l}

\define\Ind{\text{\rm Ind}}

\define\tr{\text{\rm tr}}

\define\supp{\text{\rm supp}}

\define\di{\diamond}

\define\a{\alpha}
\redefine\b{\beta}
\redefine\c{\chi}

\redefine\d{\delta}
\define\e{\epsilon}
\define\et{\eta}

\define\p{\pi}
\define\ph{\phi}

\define\r{\rho}
\define\s{\sigma}
\redefine\t{\tau}
\define\th{\theta}

\define\x{\xi}

\redefine\D{\Delta}

\define\Th{\Theta}

\define\Ph{\Phi}
\define\Ps{\Psi}

\redefine\aa{\bold a}

\define\kk{\bold k}

\define\DD{\bold D}

\define\FF{\bold F}

\define\II{\bold I}

\define\NN{\bold N}

\define\QQ{\bold Q}

\define\SS{\bold S}

\define\WW{\bold W}
\define\ZZ{\bold Z}

\define\ca{\Cal A}
\define\cb{\Cal B}
\define\cc{\Cal C}

\define\cf{\Cal F}
\define\cg{\Cal G}
\define\ch{\Cal H}

\define\cm{\Cal M}

\define\co{\Cal O}
\define\cp{\Cal P}

\define\cs{\Cal S}

\define\cu{\Cal U}

\define\cx{\Cal X}

\define\fA{\frak A}

\define\fH{\frak H}

\define\fK{\frak K}

\define\sha{\sharp}

\define\VO{AL}
\define\CHE{C}
\define\DL{DL}
\define\IRS{IRS}
\define\IW{I}
\define\KL{KL1}
\define\KLL{KL2}
\define\KO{Ko}
\define\OR{L1}
\define\QG{L2}
\define\RA{L3}
\define\UN{L4}

\head Introduction and statement of results\endhead
\subhead 0.1\endsubhead
Let $W$ be a Weyl group with standard set of generators $S$; let $\le$ be the Bruhat order
on $W$. In \cite{\KL},\cite{\KLL}, certain polynomials $P_{y,w}=\sum_{i\ge0}P_{y,w;i}u^i$ ($P_{y,w;i}\in\NN$, $u$
is an indeterminate) were defined and computed in terms of an algorithm for any $y\le w$ in $W$. These polynomials
are of interest for the representation theory of complex reductive groups, see \cite{\KL}. Let 
$\II=\{w\in W;w^2=1\}$ be the set of involutions in $W$. In this 
paper we introduce some new polynomials $P^\s_{y,w}=\sum_{i\ge0}P^\s_{y,w;i}u^i$ ($P^\s_{y,w;i}\in\ZZ$) for any
pair $y\le w$ of elements of $\II$. These new polynomials are of interest in the theory of unitary 
representations of complex reductive groups, see \cite{\VO}; they are again computable in terms of an algorithm, 
see 4.5. For $y\le w$ in $\II$ and $i\in\NN$ there is the following relation between $P_{y,w;i}$ and 
$P^\s_{y,w;i}$: there exist $a_i,b_i\in\NN$ such that $P_{y,w;i}=a_i+b_i$, $P^\s_{y,w;i}=a_i-b_i$. 

Let $\ca=\ZZ[u,u\i]$ and let $\fH$ be the free $\ca$-module with basis $(T_w)_{w\in W}$ with the unique
$\ca$-algebra structure with unit $T_1$ such that 

(i) $T_wT_{w'}=T_{ww'}$ if $l(ww')=l(w)+l(w')$ 
\nl
($l:W@>>>\NN$ is the standard length function) and $(T_s+1)(T_s-u)=0$ for all $s\in S$. Let $\fH'$ be the 
$\ca$-algebra with the same underlying $\ca$-module as $\fH$ but with multiplication defined by the rules (i) and

(ii) $(T_s+1)(T_s-u^2)=0$ for all $s\in S$.
\nl
In the course of defining the polynomials $P_{y,w}$ in \cite{\KL} a special role was played by the triple 
$(\fH\ot\fH^{opp},\fH,\,\bar{}\,:\fH@>>>\fH)$ where the middle $\fH$ is viewed as a $\fH\ot\fH^{opp}$-module via 
left and right multiplication ($\fH^{opp}$ is the algebra opposed to $\fH$) and $\,\bar{}\,:\fH@>>>\fH$ is a 
certain ring involution. To define the new polynomials $P^\s_{y,w}$ we shall instead need a triple 
$$(\fH',M,\,\bar{}\,:M@>>>M)$$
where $M$ is the free $\ca$-module with basis $(a_w)_{w\in\II}$ with a certain $\fH'$-module structure and 
$\,\bar{}\,:M@>>>M$ is a certain $\ZZ$-linear involution which are described in the following theorem (here 
$\,\bar{}\,:\fH'@>>>\fH'$ is the ring involution such that $\ov{u^nT_w}=u^{-n}T_{w\i}\i$ for all 
$w\in W,n\in\ZZ$). 

\proclaim{Theorem 0.2}(a) Consider for any $s\in S$ the $\ca$-linear map $T_s:M@>>>M$ given by

(i) $(T_s+1)(a_w)=(u+1)(a_w+a_{sw})$ if $w\in\II,sw=ws>w$;

(ii) $(T_s+1)(a_w)=(u^2-u)(a_w+a_{sw})$ if $w\in\II,sw=ws<w$;

(iii) $(T_s+1)(a_w)=a_w+a_{sws}$ if $w\in\II,sw\ne ws>w$;

(iv) $(T_s+1)(a_w)=u^2(a_w+a_{sws})$ if $w\in\II,sw\ne ws<w$.
\nl
(Note that in (iii) we have automatically $sw>w$ and in (iv) we have automatically $sw<w$.) The maps $T_s$ 
($s\in S$) define an $\fH'$-module structure on $M$.

(b) There exists a unique $\ZZ$-linear map $\,\bar{}\,:M@>>>M$ such that $\ov{u^nm}=u^{-n}\ov{m}$ for all 
$m\in M$,
$n\in\ZZ$, $\ov{a_1}=a_1$ and $\ov{(T_s+1)m}=u^{-2}(T_s+1)\ov{m}$ for all $m\in M,s\in S$. For any $w\in\II$ we 
have $\ov{a_w}=\sum_{y\in\II;y\le w}r_{y,w}a_y$ where $r_{y,w}\in\ca$ and $r_{w,w}=u^{-l(w)}$. For any $h\in\fH'$
and $m\in M$ we have $\ov{hm}=\ov{h}\ov{m}$. For any $m\in M$ we have $\ov{\ov{m}}=m$. 
\endproclaim
The proof of (a) is given in 1.8. The proof of (b), given in 2.9, is based on a sheaf theoretic construction of 
$M$, some elements of which are inspired by the geometric construction of the plus part of a universal quantized 
enveloping algebra of nonsimplylaced type given in \cite{\QG, Ch.12}.) 

Let $\uca=\ZZ[v,v\i]$ where $v$ is an indeterminate. We view $\ca$ as a subring of $\uca$ by setting $u=v^2$. Let
$\uM=\uca\ot_\ca M$. We can view $M$ as a $\ca$-submodule of $\uM$. We extend $\,\bar{}\,:M@>>>M$ to a 
$\ZZ$-linear 
map $\,\bar{}\,:\uM@>>>\uM$ in such a way that $\ov{v^nm}=v^{-n}\ov{m}$ for $m\in\uM,n\in\ZZ$. 
Let $\ufH=\uca\ot_\ca\fH$, $\ufH'=\uca\ot_\ca\fH'$. These are naturally $\uca$-algebras containing
$\fH,\fH'$ as $\ca$-subalgebras. Note that the $\fH'$-module structure on $M$ extends by $\uca$-linearity to
an $\ufH'$-module structure on $\uM$. We have the following result.

\proclaim{Theorem 0.3}(a) For any $w\in\II$ there is a unique element 
$$A_w=v^{-l(w)}\sum_{y\in\II;y\le w}P^\s_{y,w}a_y\in\uM$$
($P^\s_{y,w}\in\ZZ[u]$) such that $\ov{A_w}=A_w$, 
$P^\s_{w,w}=1$ and for any $y\in\II$, $y<w$, we have $\deg P^\s_{y,w}\le(l(w)-l(y)-1)/2$. 

(b) The elements $A_w$ ($w\in\II$) form an $\uca$-basis of $\uM$.
\endproclaim
The proof is given in 3.1, 3.2. In 3.3 we give an interpretation of $P_{y,w}^\s$ in terms of intersection 
cohomology. 

\subhead 0.4\endsubhead
For any $z\in\QQ-\{0\}$ let $M_z=\QQ\ot_\ca M$, $\fH'_z=\QQ\ot_\ca\fH'$ where $\QQ$ is viewed as an $\ca$-algebra
under $u\m z$. For $w\in W$ we write $T_w\in\fH'_z$ instead of $1\ot T_w$; for $w\in\II$ we write $a_w\in M_z$ 
instead of $1\ot a_w$. Note that $\fH'_1$ can be identified with $\QQ[W]$, the group algebra of $W$, so that for 
$w\in W$, $T_w$ becomes $w$. Now specializing 0.2(a) with $u=1$ we see that $M_1$ is a $W$-module such that

$s(a_w)=a_w+2a_{sw}$ if $w\in\II,sw=ws>w$;

$s(a_w)=-a_w$ if $w\in\II,sw=ws<w$;

$s(a_w)=a_{sws}$ if $w\in\II,sw\ne ws$.
\nl
In \S6 it is shown that the $W$-module $M_1$ is isomorphic to a direct sum of representations of $W$ induced from
one-dimensional representations of centralizers of involutions. The last direct sum has been studied in detail
by Kottwitz \cite{\KO}; in 6.4 we reformulate Kottwitz's results in terms of unipotent representations.

\subhead 0.5\endsubhead
If $X$ is a set and $f:X@>>>X$ is a map we write $X^f=\{x\in X;f(x)=x\}$. If $X$ is a finite set we write $|X|$ 
for the cardinal of $X$.

\head Contents\endhead
1. Proof of Theorem 0.2(a).

2. Proof of Theorem 0.2(b).

3. Proof of Theorem 0.3.

4. The action of $u\i(T_s+1)$ in the basis $(A_w)$.

5. Relation with two-sided cells.

6. The $W$-module $M_1$.

7. Some extensions.

\head 1. Proof of Theorem 0.2(a)\endhead
\subhead 1.1\endsubhead
Let $\kk$ be an algebraic closure of the field $\FF_p$ with $p$ elements. ($p$ is a prime number.) Let $G$ be a
connected semisimple simply connected algebraic group over $\kk$. Let $\cb$ be the variety of Borel subgroups of $G$. Then $G$ acts
on $\cb\T\cb$ by simultaneous conjugation and the set of orbits can be viewed naturally as a Coxeter group (the
Weyl group of $G$); we shall assume that this Coxeter group is $W$ of 0.1 with its standard set of generators. 
For $w\in W$ we write $\co_w$ for the corresponding $G$-orbit in $\cb\T\cb$. Let $\ph:G@>>>G$ be the Frobenius 
map for a split $\FF_p$-structure on $G$.

Let $s\in\ZZ_{>0}$ and let $q=p^s$. Then $\ph':=\ph^s:G@>>>G$ is the Frobenius map for a split $\FF_q$-structure 
on $G$ (we denote by $\FF_q$ the subfield of $\kk$ of cardinal $q$). For any $s\in S$ the $\ca$-linear map 
$T_s:M@>>>M$ defined in 0.2(i)-(iv) induces a $\QQ$-linear map $M_q@>>>M_q$ denoted again by $T_s$; it is given 
by 0.2(i)-(iv) with $u$ replaced by $q$. 

Consider the $\FF_q$-rational structure on $\cb\T\cb$ with Frobenius map 
$$F:(B,B')\m(\ph'(B'),\ph'(B)).$$
We have 
$$(\cb\T\cb)^F=\{(B,B')\in\cb\T\cb;B=\ph'(B'),B'=\ph'(B)\}.$$
If $(B,B')\in\co_w$ is fixed by $F$ then
$(\ph'(B'),\ph'(B))\in\co_w\cap\co_{w\i}$ hence $w\in\II$. On the other hand if $(B,B')\in\co_w$, $w\in\II$ then
$(\ph'(B'),\ph'(B))\in\co_w$ so that $\co_w$ is $F$-stable. We see that $(\cb\T\cb)^F=\sqc_{w\in\II}\co_w^F$.
Now if $w\in\II$, then $G$ acts transitively on $\co_w$ and this action is compatible with the $\FF_q$-structure 
on $\co_w$ given by $F$ and with the $\FF_q$-structure on $G$ given by $\ph':G@>>>G$. Hence, using Lang's 
theorem, we see that $\co_w^F\ne\em$ and that
the induced action of $G^{\ph'}$ on $\co_w^F$ is transitive (here we use
also that the stabilizer in $G$ of a point in $\co_w$ is connected). We see that the $G^{\ph'}$-orbits in
$(\cb\T\cb)^F$ are exactly the sets $\co_w^F$ with $w\in\II$. Let $\cf_q$ be the vector space of functions 
$(\cb\T\cb)^F@>>>\QQ$ which are constant on the orbits of $G^{\ph'}$. Clearly we can identify $M_q$ with $\cf_q$ 
in such a way that for $w\in\II$, $a_w$ becomes the function which is $1$ on $\co_w^F$ and is $0$ on $\co_{w'}^F$
for $w'\in\II,w'\ne w$.

Next we consider the $\FF_{q^2}$-rational structure on $\cb\T\cb$ with Frobenius map 
$(B,B')\m(\ph'{}^2(B),\ph'{}^2(B'))$ denoted again by $\ph'{}^2$. We have clearly 
$(\cb\T\cb)^{\ph'{}^2}=\sqc_{w\in W}\co_w^{\ph'{}^2}$ and this is exactly the decomposition of 
$(\cb\T\cb)^{\ph'{}^2}$ into $G^{\ph'{}^2}$-orbits. Let $\cf'_q$ be the vector space of functions 
$(\cb\T\cb)^{\ph'{}^2}@>>>\QQ$ which are constant on the orbits of $G^{\ph'{}^2}$. We define an (associative) 
algebra structure on $\cf'_q$ by $h,h'\m h*h'$ where 
$$(h*h')(B_1,B_2)=\sum_{\b\in\cb^{\ph'{}^2}}h(B_1,\b)h'(\b,B_2).$$
Clearly we can identify $\fH'_q$ with $\cf'_q$ as vector spaces in such a way that for $w\in W$, $T_w$ becomes 
the function which is $1$ on $\co_w^{\ph'{}^2}$ and is $0$ on $\co_{w'}^{\ph'{}^2}$ for $w\in W-\{w\}$. By
Iwahori \cite{\IW}, this is identification respects the algebra structures on $\cf'_q,\fH'_q$.

For $h\in\cf'_q,m\in\cf_q$ we define $h*m\in\cf_q$ by 
$$(h*m)(B_1,B_2)=\sum_{\b\in\cb^{\ph'{}^2}}h(B_1,\b)m(\b,\ph'(\b)).$$
(It may look strange that $B_2$ does not appear in the right hand side; but in fact it appears through $B_1$ since
$B_2=\ph'(B_1)$.) If $h,h'\in\cf'_q$, $m\in\cf_q$ we have
$$\align&((h*h')*m)(B_1,B_2)=\sum_{\b\in\cb^{\ph'{}^2}}(h*h')(B_1,\b)m(\b,\ph'(\b))\\&=
\sum_{\b,\b'\in\cb^{\ph'{}^2}}h(B_1,\b')h'(\b',\b)m(\b,\ph'(\b))\\&
=\sum_{\b'\in\cb^{\ph'{}^2}}h(B_1,\b')(h'*m)(\b',\ph'(\b'))=(h*(h'*m))(B_1,B_2).\endalign$$
Thus $(h*h')*m=h*(h'*m)$ so that $h,m\m h*m$ defines a $\cf'_q$-module structure on $\cf_q$. (Note that 
$T_1*m=m$ for $m\in\cf_q$.) 

\subhead 1.2\endsubhead
Let $s\in S$, $w\in\II$. We have $T_s*a_w=\sum_{w'\in\II}N_{s,w,w'}a_{w'}$, where
$$\align&N_{s,w,w'}=|\{\b\in\cb^{\ph'{}^2};(\b,\ph'(\b))\in\co_w,(C,\b)\in\co_s\}|\\&
=|\{(\b,\b')\in\co_w^F;(C,\b)\in\co_s,(\b',C')\in\co_s\}|\endalign$$
for any $(C,C')\in\co_{w'}^F$. To simplify notation we write $N_{w'}$ instead of $N_{s,w,w'}$ since $s,w$ are fixed.
In 1.3-1.6 we compute the number $N_{w'}$ for a fixed $(C,C')\in\co_{w'}^F$, $w'\in\II$.

\subhead 1.3\endsubhead
In this subsection we assume that $sw\ne ws>w$. Using \cite{\DL, 1.6.4} we see that $l(sws)\ne l(w)$ hence 
$l(sws)=l(w)+2$. If $N_{w'}\ne0$ then there exists $(\b,\b')\in\co_w$ such that 
$(C,\b)\in\co_s,(\b',C')\in\co_s$. Hence $w'=sws$. Conversely, assume that $w'=sws$. There is a unique 
$(\b,\b')\in\co_w$ such that $(C,\b)\in\co_s,(\b',C')\in\co_s$.
We have $(\ph'(\b'),\ph'(\b))\in\co_w$, $(C,\ph'(\b'))\in\co_s,(\ph'(\b),C')\in\co_s$. (We use that $\ph'(C)=C'$,
$\ph'(C)=C'$.) By uniqueness we must have $\ph'(\b)=\b'$, $\ph'(\b')=\b$. Thus $(\b,\b')\in\co_w^F$. We see that
$N_{sws}=1$ so that $T_s*a_w=a_{sws}$.

\subhead 1.4\endsubhead
In this subsection we assume that $sw=ws>w$. If $N_{w'}\ne0$ then there exists  $(\b,\b')\in\co_w$ such that 
$(C,\b)\in\co_s,(\b',C')\in\co_s$. Hence $(C,\b')\in\co_{sw}$ (we use that $l(sw)>l(w)$) and 
$(C,C')\in\co_{sw}\cup\co_w$ (we use that $sws=w$); hence $w'=sw$ or $w'=w$.

Assume first that $w'=w$. We set 
$$Z=\{(\b,\b')\in\co_w;(C,\b)\in\co_s,(\b',C')\in\co_s\}.$$
We claim that the first projection $Z@>>>Y:=\{\b\in\cb;(C,\b)\in\co_s\}$ is an isomorphism so that $Z$ is an 
affine line. Let $\b\in Y$. It is enough to show that there is a unique $\b'\in\cb$ such that
$(\b,\b')\in\co_w,(\b',C')\in\co_s$. Since $l(ws)>l(w)$ it is enough to show that $(\b,C')\in\co_{ws}$; but from
$(\b,C)\in\co_s, (C,C')\in\co_w$, $l(sw)>l(w)$ we do indeed deduce that $(\b,C')\in\co_{ws}$, as desired. This 
proves our claim. Now the restriction of $F$ to $Z$ is an $\FF_q$-rational structure on an affine line hence it 
has exactly $q$ fixed points. It follows that $N_w=q$.

Assume next that $w'=sw=sw$. We set 
$$Z'=\{(\b,\b')\in\co_w;(C,\b)\in\co_s\cup\co_1,(\b',C')\in\co_s\cup\co_1\}.$$
We claim that the first projection $Z'@>>>Y':=\{\b\in\cb;(C,\b)\in\co_s\cup\co_1\}$ is an isomorphism so that 
$Z'$ is a projective line. It is enough to show that for any $\b\in Y'$ the set 
$$\Xi_\b=\{\b'\in\cb;(\b,\b')\in\co_w,(\b',C')\in\co_s\cup\co_1\}$$ 
has exactly one element. Let $C_1$ be the unique element of $\cb$ such that $(C,C_1)\in\co_s,(C_1,C')\in\co_w$ 
(we use that $(C,C')\in\co_{sw}$ and $l(sw)>l(w)$). Note that $C_1\in Y'$. If $\b=C_1$ then $C'\in\Xi_\b$; if 
$\b'\in\Xi_b-\{C'\}$ then $(\b',C')\in\co_s$, $(C',C_1)\in\co_w$ hence $(\b',C_1)\in\co_{sw}$ (since $l(sw)>l(w)$)
contradicting $(\b',\b)\in\co_w$. Thus $|\Xi_{C_1}|=|\{C'\}|=1$. Now assume that $\b\in Y'-\{C_1\}$. From 
$\b\ne C_1$, $(C,\b)\in\co_s\cup\co_1$, $(C,C')\in\co_{sw}$ we deduce that $(\b,C')\in\co_{sw}$. Hence there is a
unique $\b'_0\in\cb$ such that $(\b,\b'_0)\in\co_w$, $(\b'_0,C')\in\co_s$. We have $\b'_0\in\Xi_\b$. Conversely,
if $\b'\in\Xi_\b$ we have $\b'\ne C'$ (since $(\b,C')\in\co_{sw}$); thus $\b'=\b'_0$. We see that
$|\Xi_\b|=|\{b'_0\}|=1$. This proves our claim.

Now the restriction of $F$ to $Z'$ is an $\FF_q$-rational structure on a projective line hence it has exactly 
$q+1$ fixed points. If $(\b,\b')\in Z^F$ then we have necessarily $C\ne\b$ and $C'\ne\b'$. Indeed, if $C=\b$ then
$\ph'(C)=\ph'(\b)$ hence $C'=\b'$ and $(\b,\b')=(C,C')\in\co_{sw}$ contradicting $(\b,\b')\in\co_w$. Similarly, 
if $C'=\b'$ then $\ph'(C')=\ph'(\b')$ hence $C=\b$ and $(\b,\b')=(C,C')\in\co_{sw}$ contradicting 
$(\b,\b')\in\co_w$. Thus 
$$q+1=|Z'{}^F|=|\{(\b,\b')\in\co_w^F;(C,\b)\in\co_s,(\b',C')\in\co_s\}|=N_{sw}.$$
We see that $T_s*a_w=qa_w+(q+1)a_{sw}$.

\subhead 1.5\endsubhead
In this subsection we assume that $sw\ne ws<w$. Using \cite{\DL, 1.6.4} we see that $l(sws)\ne l(w)$ hence 
$l(sws)=l(w)-2$. We have $s(sws)\ne(sws)s>sws$. Applying 1.3 to $sws$ instead of $w$ we obtain $T_s*a_{sws}=a_w$.
Applying $T_s$ to the last equality and using the equation $T_s^2=(q^2-1)T_s+q^2:\cf_q@>>>\cf_q$, we obtain 

$(q^2-1)T_s*a_{sws}+q^2a_{sws}=T_s*a_w$
\nl
that is, $T_s*a_w=(q^2-1)a_w+q^2a_{sws}$.

\subhead 1.6\endsubhead
In this subsection we assume that $sw=ws<w$. We have $s(sw)=(sw)s>sw$. Applying 1.4 to $sw$ instead of $w$ we 
obtain $T_s*a_{sw}=qa_{sw}+(q+1)a_w$. Applying $T_s$ to the last equality and using the equation 
$T_s^2=(q^2-1)T_s+q^2:\cf_q@>>>\cf_q$ we obtain $(q^2-1)T_s*a_{sw}+q^2a_{sw}=qT_s*a_{sw}+(q+1)T_s*a_w$ that is
$$\align&(q+1)T_s*a_w=(q^2-1-q)T_s*a_{sw}+q^2a_{sw}\\&
=(q^2-q-1)qa_{sw}+(q^2-q-1)(q+1)a_w+q^2a_{sw}.\endalign$$
Dividing by $q+1$ we obtain $T_s*a_w=(q^2-q-1)a_w+(q^2-q)a_{sw}$.

\subhead 1.7\endsubhead
From the results in 1.3-1.6 we see that the operators $T_s:M_q@>>>M_q$ ($s\in S$) given by the formulas 
0.2(i)-(iv) with $u$ replaced by $q$ define an $\fH'_q$-module structure on $M_q$.

\subhead 1.8\endsubhead
We now prove 0.2(a). We shall use the following obvious fact.

(a) If $m\in M$ has zero image in $M_{p^s}$ for $s=1,2,\do$ then $m=0$.
\nl
Let $s,t\in S,s\ne t$ and let $k$ be the order of $st$ in $W$. Let $m\in M$. Let 
$m'=(T_sT_tT_s\do)m-(T_tT_sT_t\do)m\in M$ (both products have $k$ factors). From 1.7 we see that $m'$ has zero 
image in $M_{p^s}$ for $s=1,2,\do$; hence by (a) we have $m'=0$.

Now let $s\in S$ and let $m\in M$. Let $m'=T_s^2m-(u^2-1)T_sm-u^2m\in M$. From 1.7 we see that $m'$ has zero 
image in $M_{p^s}$ for $s=1,2,\do$; hence by (a) we have $m'=0$.

We see that 0.2(a) holds.

\head 2. Proof of Theorem 0.2(b)\endhead
\subhead 2.1\endsubhead
We preserve the notation of 1.1. We fix a prime number $l\ne p$. For any complex $\fK$ of constructible 
$\ql$-sheaves on an algebraic variety we denote by $\ch^i\fK$ the $i$-th cohomology sheaf of $\fK$ and by $D\fK$
the Verdier dual of $\fK$.

Let $\cc_0$ be the category whose objects are the constructible $G$-equivariant $\ql$-sheaves on $\cb\T\cb$; the 
morphisms in $\cc_0$ are morphisms of $G$-equivariant $\bbq$-sheaves. Let $Vec$ be the category of finite-dimensional 
$\ql$-vector spaces. If $\cs\in\cc_0$ and $x\in\cb\T\cb$, let $\cs_x$ be the stalk of $\cs$ at $x$; 
for any $w\in W$ there is a well defined object $V^\cs_w\in Vec$ which is canonically isomorphic to $\cs_x$ for 
any $x\in\co_w$. Note that $\cs\m V^\cs_w$ is a functor $\cc_0@>>>Vec$. For $w\in W$ let $\SS_w$ be the object of
$\cc_0$ which is $\ql$ on $\co_w$ and is $0$ on $\cb\T\cb-\co_w$. For $w\in W$ let $\SS_w^\sha$ be the 
intersection cohomology complex of the closure $\bco_w$ of $\co_w$ with coefficients in $\ql$ (on $\co_w$) 
extended by zero on $\cb\T\cb-\bco_w$. We have $\ch^i\SS_w^\sha\in\cc_0$ ($i\in\ZZ$); moreover $\ch^i\SS_w^\sha$ 
is zero for large $i$ and is zero unless $i\in2\NN$, see \cite{\KLL, Thm. 4.2}. If $\cs\in\cc_0$, then $\ch^iD\cs\in\cc_0$ 
($i\in\ZZ$); moreover $\ch^iD\cs$ is zero for all but finitely many $i$.

Let $\cc_1$ be the category whose objects are pairs $(\cs,\Ps)$ where $\cs\in\cc_0$ and $\Ps$ is an isomorphism 
$\ph^{2*}\cs@>\si>>\cs$ in $\cc_0$. A morphism between two objects $(\cs,\Ps),(\cs',\Ps')$ of $\cc_1$ is a 
morphism $\cs@>>>\cs'$ in $\cs_0$ which is compatible with $\Ps,\Ps'$. If $(\cs,\Ps)\in\cc_1$ then for any 
$w\in W$, $\Ps$ induces a linear isomorphism $\Ps_w:V^\cs_w@>>>V^\cs_w$ (Note that 
$V^{\ph^{2*}\cs}_w=V^{\cs}_w$ since $\co_w$ is $\ph$-stable.) If $z\in\ql^*$ and $(\cs,\Ps)\in\cc_1$ then 
$(\cs,z\Ps)\in\cc_1$. If $(\cs,\Ps)\in\cc_1$ then for any $i$ let $\Ps^{(i)}:\ph^{2*}\ch^iD\cs@>\si>>\ch^iD\cs$ 
be the inverse of the isomorphism $\ch^iD\cs@>>>\ph^{2*}\ch^iD\cs$ induced by $D$; thus 
$(\ch^iD\cs,\Ps^{(i)})\in\cc_1$. 

For any $w\in W$ and any $z\in\ql^*$ we have $(\SS_w,t_z)\in\cc_1$ where $t_z:\ph^{2*}\SS_w@>\si>>\SS_w$ is 
multiplication by $z$ (we have $\ph^{2*}\SS_w=\SS_w$).

Let $\cc_2$ be the full category of $\cc_1$ whose objects are the pairs $(\cs,\Ps)\in\cc_1$ such that for any 
$w\in W$ all eigenvalues of $\Ps_w:V^\cs_w@>>>V^\cs_w$ are even powers of $p$. For example if $w\in W,k\in\ZZ$ 
then $(\SS_w,t_{p^{2k}})\in\cc_2$.

\proclaim{Lemma 2.2} For any $(\cs,\Ps)\in\cc_2$ and any $i\in\ZZ$ we have $(\ch^iD\cs,\Ps^{(i)})\in\cc_2$.
\endproclaim
It suffices to show that $(\ch^iD\SS_w,t_1^{(i)})\in\cc_2$ for any $w\in W,i\in\ZZ$. We can assume that $w\in W$ 
is fixed and that the statement in the previous sentence holds when $w$ is replaced by any $w'\in W$, $w'<w$. Let
$j^w:\ph^{2*}\SS_w^\sha@>\si>>\SS_w^\sha$ be the unique isomorphism such that the
induced isomorphism $V^{\ch^0\SS_w^\sha}_w@>>>V^{\ch^0\SS_w^\sha}_w$ is the identity map.
By \cite{\KLL, Thm. 4.2}, $j^w$ induces for any $h\in\ZZ$ an isomorphism 
$j^w_h:\ph^{2*}\ch^h\SS_w^\sha@>\si>>\ch^h\SS_w^\sha$ such that 

(a) $(\ch^h\SS_w^\sha,j^w_h)\in\cc_2$.
\nl
Let $K$ be the restriction of $\SS_w^\sha$ to $\bco_w-\co_w$ extended by $0$ on the complement of $\bco_w-\co_w$.
Now $j^w$ induces an isomorphism $d:\ph^{2*}K@>\si>>K$. Let $d':\ph^{2*}DK@>\si>>DK$ be the inverse of the
isomorphism $DK@>>>\ph^{2*}DK$ induced by $d$; this induces isomorphisms $d'_h:\ph^{2*}\ch^hDK@>\si>>\ch^hDK$ for
$h\in\ZZ$. Since $\supp(K)\sub\bco_w-\co_w$ we see using (a) and the induction hypothesis that 

(b) $(\ch^hDK,d'_h)\in\cc_2$ for $h\in\ZZ$.
\nl
Now $j^w$ induces an isomorphism $D\SS_w^\sha@>>>\ph^{2*}D\SS_w^\sha$ whose inverse is an isomorphism
$j':\ph^{2*}D\SS_w^\sha@>\si>>D\SS_w^\sha$. This induces for any $h\in\ZZ$ an isomorphism
$j'_h:\ph^{2*}\ch^hD\SS_w^\sha@>\si>>\ch^hD\SS_w^\sha$. We can identify $D\SS_w^\sha= \SS_w^\sha[2m]$ for some 
$m\in\ZZ$ in such a way that $j'_h$ becomes $p^{2m'}j^w_{h+2m}$ for some $m'\in\ZZ$. Using (a) we deduce that

(c) $(\ch^hD\SS_w^\sha,j'_h)\in\cc_2$.
\nl
Using (c),(b) and the long exact sequence of cohomology sheaves associated to the exact triangle consisting of 
$D\SS_w,\D\SS_w^\sha,DK$ (which is obtained from the exact triangle consisting of $K,\SS_w^\sha,\SS_w$) we deduce
that $(\ch^hD\SS_w,t_1^{(h)})\in\cc_2$ for any $h\in\ZZ$. This completes the inductive proof.

\subhead 2.3\endsubhead
Define $\s:\cb\T\cb@>>>\cb\T\cb$ by $(B,B')\m(B',B)$. Then 
$$\tph:=\ph\s=\s\ph:\cb\T\cb@>>>\cb\T\cb$$ 
is the Frobenius map for an $\FF_p$-rational structure on $\cb\T\cb$ such that $\tph^2=\ph^2$.

Let $\cc$ be the category whose objects are pairs $(\cs,\Ph)$ where $\cs\in\cc_0$ and $\Ph$ is an isomorphism 
$\tph^*\cs@>>>\cs$ in $\cc$ such that, setting $\Ps=\tph^*(\Ph)\Ph:\tph^{2*}\cs@>>>\cs$, we have 
$(\cs,\Ps)\in\cc_2$ (note that $\tph^{2*}\cs=\ph^{2*}\cs$). A morphism between two objects 
$(\cs,\Ph),(\cs',\Ph')$ of $\cc$ is a 
morphism $\cs@>>>\cs'$ in $\cs_0$ which is compatible with $\Ph,\Ph'$. Note that if $(\cs,\Ph)\in\cc$, then 
$(\cs,-\Ph)\in\cc$. For $(\cs,\Ph)\in\cc$ and $w\in\II$, $\Ph$ induces a linear isomorphism 
$\Ph_w:V^\cs_w@>>>V^\cs_w$. (Note that $V^{\tph^*\cs}_w=V^\cs_w$ since $\co_w$ is $\tph$-stable.)

For any $(\cs,\Ps)\in\cc_2$ let $\Ph:\tph^*(\cs\op\tph^*\cs)@>>>\cs\op\tph^*\cs$ (that is 
$\Ph:\tph^*\cs\op\tph^{2*}\cs@>>>\cs\op\tph^*\cs$) be the isomorphism whose restriction to
$\tph^*\cs$ is $0\op 1$ and whose restriction to $\tph^{2*}\cs$ is $\Ps\op 0$. We set
$\Th(\cs,\Ps)=(\cs\op\tph^*\cs,\Ph)$. Note that $\Th(\cs,\Ps)\in\cc$. 

Let $E$ be the subset of $\ql$ consisting of $p^n,-p^n$ ($n\in\ZZ$). For $w\in\II$, $z\in E$ let
$\t_z:\tph^*\SS_w@>\si>>\SS_w$ be the isomorphism such that $\Ph_w:V^{\SS_w}_w@>>>V^{\SS_w}_w$ is multiplication 
by $z$. We have $(\SS_w,\t_z)\in\cc$.

Let $K(\cc)$ be the Grothendieck group of $\cc$. We have the following result.

\proclaim{Lemma 2.4} If $(\cs,\Ph)\in\cc$ then in $K(\cc)$, $(\cs,\Ph)$ is a $\ZZ$-linear combination of elements
$(\SS_w,\t_z)$ ($w\in\II,z\in E$) and of elements $\Th(\cs,\Ps)$ for various $(\cs,\Ps)\in\cc_2$. 
\endproclaim
We can assume that the set $J:=\{w\in W;\cs|_{\co_w}\ne0\}$ consists either of (i) a single element of $\II$ or 
(ii) of two distinct elements $w',w''$ whose product is $1$ and that for any $w\in J$ we have $\cs|_{\co_w}=\ql$. 
In case (i) we have $(\co,\Ph)=(\SS_w,\t_z)$ where $J=\{w\}$, $z\in E$; in case (ii) we have 
$(\co,\Ph)\cong\Th(\SS_w,t_{p^{2n}})$ where $w\in J$ and $n\in\ZZ$. The lemma is proved.

\subhead 2.5\endsubhead
Let $K'(\cc)$ be the subgroup of $K(\cc)$ generated by the elements of the form $\Th(\cs,\Ps)$ for various 
$(\cs,\Ps)\in\cc_2$ and by the elements of the form $(\cs,\Ph)+(\cs,-\Ph)$ with $(\cs,\Ph)\in\cc$. Let 
$\bK(\cc)=K(\cc)/K'(\cc)$. From 2.4 we see that 

(a) {\it the abelian group $\bK(\cc)$ is generated by the elements $(\SS_w,\t_{p^n})$ ($w\in\II,n\in\ZZ$).}
\nl
We regard $K(\cc)$ as an $\ca$-module where $u^n(\cs,\Ph)=(\cs,p^n\Ph)$ for $n\in\ZZ$. Then $K'(\cc)$ is an 
$\ca$-submodule of $K(\cc)$ hence $\bK(\cc)$ inherits an $\ca$-module structure from $K(\cc)$.

\subhead 2.6\endsubhead
Let $s$ be an odd integer $\ge1$ and let $q=p^s$. Let 
$$F=\ph^s\s=\s\ph^s=\tph^s:\cb\T\cb@>>>\cb\T\cb$$ 
be as in 1.1. For $(\cs,\Ph)\in\cc$ we define a function $\c_s(\cs,\Ph)\in\cf_q$ (see 1.1) as follows. For 
$x\in(\cb\T\cb)^F$, $\c_s(\cs,\Ph)(x)$ is the trace of the composition 
$$\cs_x=\cs_{\tph^s(x)}@>\Ph>>\cs_{\tph^{s-1}(x)}@>\Ph>>\do@>\Ph>>\cs_{\tph(x)}@>\Ph>>\cs_x$$
or equivalently, the trace of $\Ph_w^s:V^\cs_w@>>>V^\cs_w$ where $w\in\II$ is such that $x\in\co_w^{\tph^s}$.
Clearly, $(\cs,\Ph)\m\c_s(\cs,\Ph)$ defines a group homomorphism $\c_s:K(\cc)@>>>\cf_q$ such that 
$\c_s(u\x)=p^s\c_s(\x)$ for all $\x\in K(\cc)$ and such that $\c_s(\SS_w,\t_1)=a_w$ for all $w\in\II$.

If $(\cs,\Ps)\in\cc_2$, the function $\c_s(\Th(\cs,\Ph)):(\cb\T\cb)^F@>>>\bbq$ is $0$ (its value at 
$x\in(\cb\T\cb)^F$ is, from the definitions, the trace of a linear map of the form $A\op B@>>>A\op B$, 
$(a,b)\m(T(b),T'(a))$ where $T:B@>>>A,T':A@>>>B$ are linear maps). From the definition we see also that if 
$(\cs,\Ph)\in\cc$ then $\c_s(\cs,\Ph)+\c_s(\cs,-\Ph)=0$. Thus, $\c_s:K(\cc)@>>>\cf_q$ maps $K'(\cc)$ to $0$ hence
it induces a group homomorphism $\bK(\cc)@>>>\cf_q$ denoted again by $\c_s$. It satisfies $\c_s(u\x)=p^s\c_s(\x)$
for all $\x\in\bK(\cc)$. We show:

(a) {\it if $\x\in\bK(\cc)$ is such that $\c_s(\x)=0$ for all odd integers $s\ge1$ then $\x=0$.}
\nl
By 2.4 we have $\x=\sum_{w\in\II,n\in\ZZ}c_{w,n}(\SS_w,\t_{p^n})$ where $c_{w,n}\in\ZZ$ are zero for all but 
finitely many $(w,n)$. Applying $\c_s$ we obtain $0=\sum_{w\in\II,n\in\ZZ}c_{w,n}p^{ns}a_w$. Since the $a_w$ form
a basis of $\cf_{p^s}$ we deduce $\sum_{w\in\II,n\in\ZZ}c_{w,n}p^{ns}=0$ for any $w\in\II$. Since this holds for 
$s=1,3,\do$ we see that $c_{w,n}=0$ for any $w\in\II,n\in\ZZ$, proving (a).

We show:

(b) {\it The elements $(\SS_w,\t_1)$ $(w\in\II)$ form an $\ca$-basis of $\bK(\cc)$.}
\nl
The fact that they generate the $\ca$-module $\bK(\cc)$ follows from 2.4. The fact that they are linearly
independent over $\ca$ follows from the proof of (a).

\subhead 2.7\endsubhead
If $(\cs,\Ph)\in\cc$ then for any $i$ let $\Ph^{(i)}:\tph^*\ch^iD\cs@>\si>>\ch^iD\cs$ be the inverse of the 
isomorphism $\ch^iD\cs@>>>\tph^*\ch^iD\cs$ induced by $D$. Note that if $z\in E$ (see 2.3) then 
$(z\Ps)^{(i)}=z\i\Ps^{(i)}$. Moreover, setting $\Ps=\tph^*(\Ph)\Ph:\tph^{2*}\cs@>\si>>\cs$, we have 
that $\tph^*(\Ph^{(i)})\Ph^{(i)}:\tph^{2*}\cs@>>>\cs$ is equal to $\Ps^{(i)}$ (as in 2.1). Using now 2.2 we see 
that

(a) $(\ch^iD\cs,\Ph^{(i)})\in\cc$.
\nl
Clearly there is a well defined $\ZZ$-linear map $\DD:K(\cc)@>>>K(\cc)$ such that
$$\DD(\cs,\Ph)=\sum_{i\in\ZZ}(-1)^i(\ch^iD\cs,\Ph^{(i)})$$ 
for any $(\cs,\Ph)\in\cc$. If $(\cs,\Ps)\in\cc_2$ and $\Th(\cs,\Ps)=(\cs\op\tph^*\cs,\Ph)$ then for any $i$ we 
have $(\ch^iD(\cs\op\tph^*\cs),\Ph^{(i)})=\Th(\ch^iD\cs,\Ps^{(i)})$ where $(\ch^iD\cs,\Ps^{(i)})\in\cc_2$ (see 
2.2). Moreover if $(\cs,\Ph)\in\cc$ then for any $i$ we have $(\ch^iD\cs,(-\Ph)^{(i)})=(\ch^iD\cs,-\Ph^{(i)})$. 
It follows that $\DD$ carries $K'(\cc)$ into itself hence it induces a $\ZZ$-linear map $\bK(\cc)@>>>\bK(\cc)$ 
denoted again by $\DD$. Note that $\DD(u^n\xi)=u^{-n}\DD(\xi)$ for any $\xi\in\bK(\cc)$ and any $n\in\ZZ$.

Since $\co_1$ is closed, smooth, of pure dimension $\nu:=\dim\cb$, we have from the definitions
$$\DD(\SS_1,\t_1)=u^{-\nu}(\SS_1,\t_1).$$

\subhead 2.8\endsubhead
Now let $t\in S$. We have $\bco_t=\co_t\cup\co_1$. Let
$$Y=\{(B_1,B_2,B_3,B_4)\in\cb^4;(B_1,B_2)\in\bco_t,(B_3,B_4)\in\bco_t\}.$$ 
Define $\tph:Y@>>>Y$ by 
$$\tph(B_1,B_2,B_3,B_4)=(\ph(B_4),\ph(B_3),\ph(B_2),\ph(B_1)).$$
This is the Frobenius map for an $\FF_p$-rational structure on $Y$.
Define $\p,\p':Y@>>>\cb\T\cb$ by 
$$\p(B_1,B_2,B_3,B_4)=(B_2,B_3),\qua \p'(B_1,B_2,B_3,B_4)=(B_1,B_4).$$ 
We have $\p\tph=\tph\p,\p'\tph=\tph\p'$. (The $\tph$ to the left of $\p$ or $\p'$ is as in 2.3.)
For $\cs\in\cc_0$ and $i\in\ZZ$ let $\cs^i=R^i\p'_!\p^*\cs$; note that $\cs^i\in\cc_0$.
Let $(\cs,\Ps)\in\cc_2$. For $i\in\ZZ$, $\Ps:\ph^{2*}\cs@>>>\cs$ induces an isomorphism 
$\ph^{2*}\p^*\cs@>>>\p^*\cs$ (since $\p\ph^2=\ph^2\p$) and this induces for any $i$ an isomorphism 
${}^i\Ps:\ph^{2*}\cs^i@>>>\cs^i$ (since $\p'\ph^2=\ph^2\p'$). A standard argument shows that 
$(\cs^i,{}^i\Ps)\in\cc_2$. It follows that if $(\cs,\Ph)\in\cc$ and $i\in\ZZ$, then the isomorphism 
${}^i\Ph:\tph^*\cs^i@>>>\cs^i$ induced by $\Ph$ satisfies $(\cs^i,{}^i\Ph)\in\cc$. (Setting 
$\Ps=\tph^*(\Ph)\Ph:\tph^{2*}\cs@>\si>>\cs$, we have ${}^i\Ps=\tph^*({}^i\Ph)({}^i\Ph):\tph^{2*}\cs@>\si>>\cs$.)
Hence there is a well defined $\ZZ$-linear map 
$\th_t:K(\cc)@>>>K(\cc)$ such that $\th_t(\cs,\Ph)=\sum_i(-1)^i(\cs^i,{}^i\Ph)$ for all $(\cs,\Ph)\in\cc$.
From the definitions we have $\th_t(u^n\x)=u^n\th_t(\x)$ for any $\x\in K(\cc)$, $n\in\ZZ$.
From the known properties of Verdier duality we have that
$$\DD(\th_t(\x))=u^{-2}\th_t(\DD(\x))\text{ for all }\x\in K(\cc).$$
(We use that $\p'$ is proper and that $\p$ is smooth with connected fibres of dimension $2$.)

If $(\cs,\Ps)\in\cc_2$ and $\Th(\cs,\Ps)=(\cs\op\tph^*\cs,\Ph)$ then for any $i$ we have

$(\cs\op\tph^*\cs)^i,{}^i\Ph)=\Th(\cs^i,{}^i\Ps)$.
\nl
Moreover if $(\cs,\Ph)\in\cc$ then for any $i$ we have 

$(\cs^i,{}^i(-\Ph))=(\cs^i,-{}^i\Ph)$. 
\nl
It follows that $\th_t$ carries $K'(\cc)$ into itself hence it induces an
$\ca$-linear map $\bK(\cc)@>>>\bK(\cc)$ denoted again by $\th_t$.

Now let $s$ be an odd integer $\ge1$ and let $q=p^s$. We define a linear map $\th_{t,s}:\cf_q@>>>\cf_q$ by 
$f\m f'$ where
$$f'(B_1,B_4)=\sum_{(B_2,B_3)\in(\cb\T\cb)^{\tph^s};(B_1,B_2)\in\bco_t,(B_3,B_4)\in\bco_t}f(B_2,B_3)$$
for any $(B_1,B_4)\in(\cb\T\cb)^{\tph^s}$. 

Let $(\cs,\Ph)\in\cc$ and define $(\cs^i,{}^i\Ph)\in\cc$ as above. Using Grothendieck's sheaves-functions
dictionary, we see that 
$$\th_{t,s}(\c_s(\cs,\Ph))=\sum_i(-1)^i\c_s(\cs_i,{}^i\Ph).$$
It follows that
$$\th_{t,s}(\c_s(\x))=\c_s(\th_t(\x))$$
for any $\x\in\bK(\cc)$. From the definition of the $\ch'_q$-module structure on $\cf_q$ in 1.1 we see that
$\th_{t,s}(m)=(T_t+1)(m)$ for any $m\in\cf_q$. Thus for any $\x\in\bK(\cc)$ we have 
$$\c_s(\th_t(\x))=(T_t+1)(\c_s(\x)).$$

\subhead 2.9\endsubhead
Using 2.6(b) we identify $M=\bK(\cc)$ as $\ca$-modules in such a way that for any $w\in\II$, the element
$a_w\in M$ becomes the element $(\SS_w,\t_1)$ of $\bK(\cc)$.
Then $\DD$ in 2.7 and $\th_t$ in 2.8 become $\ZZ$-linear maps $M@>>>M$ (denoted again by $\DD,\th_t$).
From 2.8 we see that for any $s=1,3,5,\do$ and any $m\in M$, the image of $\th_t(m)$ in $M_{p^s}$ is
equal to the image of $(T_t+1)m$ in $M_{p^s}$. It follows that 

$\th_t(m)=(T_t+1)m$ in $M$. 
\nl
Note also that $\th_t:M@>>>M$ is $\ca$-linear while $\DD(u^nm)=u^{-n}\DD(m)$ for all $m\in M$, $n\in\ZZ$.
From 2.8 we see that $\DD(\th_t(m))=u^{-2}\th_t(\DD(m))$ for all $m\in M$. Equivalently we have
$\DD((T_t+1)(m))=u^{-2}(T_t+1)(\DD(m))$ for all $m\in M$. (Here $t$ is any element of $S$.) From 2.7
we have $\DD(a_1)=u^{-\nu}a_1$. We now define $\,\bar{}\,:M@>>>M$ by $m\m u^\nu\DD(m)$. This has the properties
described in the first sentence of 0.2(b). Thus the existence part of that sentence is established. To prove the 
uniqueness part of that sentence it is enough to verify the following statement.

{\it Let $f:M@>>>M$ be a $\ZZ$-linear map such that $f(u^nm)=u^{-n}\ov{m}$ for all $m\in M$, $n\in\ZZ$, 
$f(a_1)=0$ and $f((T_t+1)m)=u^{-2}(T_t+1)f(m)$ for all $m\in M,t\in S$. Then $f=0$.}
\nl
We must show that $f(a_w)=0$ for any $w\in\II$. We can assume that $w\ne1$ and that $f(a_{w'})=0$ for
any $w'\in\II$ such that $w'<w$. We can find $t\in S$ such that $wt<w$. If $tw\ne wt$ then applying 0.2(iii) with
$w,s$ replaced by $twt,t$ we have 
$$f(u^2a_{twt}+u^2a_w)=f((T_t+1)(a_{twt}))=u^{-2}(T_t+1)f(a_{twt})=0$$ 
(since $twt<w$) hence $u^{-2}f(a_w)=0$ so that $f(a_w)=0$. If $tw=wt$ then applying 0.2(i) with $w,s$ replaced by
$wt,t$ we have 
$$f((u+1)(a_{wt}+a_w))=f((T_t+1)(a_{wt}))=0$$ 
hence $(u\i+1)f(a_w)=0$ so that $f(a_w)=0$. This completes the proof of the first sentence in 0.2(b).
The second sentence in 0.2(b) follows from the fact that $D\SS_w$ has support contained
in $\bco_w$ and from the fact that $\dim(\co_w)=l(w)+\nu$. We prove the third sentence in 0.2(b).
Let $\fH'_0$ be the set of all $h\in\fH'$ such that $\ov{hm}=\ov{h}\ov{m}$ for any $m\in M$.
Clearly $\fH'_0$ is an $\ca$-subalgebra with $1$ of $\fH'$. By definition, $\fH'_0$ contains
$T_t+1$ for any $t\in S$. Since the elements $T_t+1$ generate $\fH'$ as an $\ca$-algebra we see that 
$\fH'_0=\fH'$, as desired. We prove the fourth sentence in 0.2(b). Define $f':M@>>>M$ by
$f'(m)=\ov{\ov{m}}$. This is an $\fH'$-linear map $M@>>>M$ such that $f'(a_1)=a_1$. We must show that
$f'=1$. It is enough to show that $f'=1$ after the scalars are extended to $\ZZ[u,u\i,(u+1)\i]$. But this follows
the fact that $M$ (with scalars thus extended) is generated by $a_1$ as an $\fH'$-module. (We use the formulas in
0.2(a).) This completes the proof of 0.2(b) hence that of Theorem 0.2.

\head 3. Proof of Theorem 0.3\endhead
\subhead 3.1\endsubhead
In this section we fix $w\in\II$. Recall that we have morphisms $\s,\ph,\tph=\s\ph=\ph\s$ of $\cb\T\cb$ into 
itself, see 2.3. Clearly there is a unique isomorphism 
$k^w:\s^*\SS_w^\sha@>\si>>\SS_w^\sha$, (resp. $m^w:\ph^*\SS_w^\sha@>\si>>\SS_w^\sha$, 
$i^w:\tph^*\SS_w^\sha@>\si>>\SS_w^\sha$)
such that the induced isomorphisms

$(k^w_h)_y:V^{\ch^h\SS_w^\sha}_y@>\si>>V^{\ch^h\SS_w^\sha}_y$ 

(resp. $(m^w_h)_y:V^{\ch^h\SS_w^\sha}_y@>\si>>V^{\ch^h\SS_w^\sha}_y$, 
$(i^w_h)_y:V^{\ch^h\SS_w^\sha}_y@>\si>>V^{\ch^h\SS_w^\sha}_y$)
\nl
defined for any $y\in\II,y\le w$ and $h\in2\NN$ satisfy $(k^w_0)_w=1$ (resp. $(m^w_0)_w=1$, $(i^w_0)_w=1$).
We have $(i^w_h)_y=(k^w_h)_y(m^w_h)_y=(m^w_h)_y(k^w_h)_y$. By \cite{\KLL, Thm. 4.2}, $(m^w_h)_y$ is equal to $p^{h/2}$ 
times a unipotent transformation. Clearly, $(k^w_h)_y$ has square equal to $1$. In particular, for any 
$h\in2\NN$, we have
$$\align&(\ch^h\SS_w^\sha,i^w_h)\in\cc,\\&
(\ch^h\SS_w^\sha,i^w_h)=\sum_{y\in\II;y\le w}P^\s_{y,w;h/2}u^{h/2}\SS_y
=\sum_{y\in\II;y\le w}P^\s_{y,w}\SS_y\in\bK(\cc)\tag a\endalign$$
where
$$P^\s_{y,w;h/2}=\tr((k^w_h)_y:V^{\ch^h\SS_w^\sha}_y@>>>V^{\ch^h\SS_w^\sha}_y)\in\ZZ,$$ 
$$P^\s_{y,w}=\sum_{h\in2\NN}P^\s_{y,w;h/2}u^{h/2}\in\ZZ[u].$$
From the definition of $\SS_w^\sha$ we have (with notation of 2.7):
$$\sum_{h\in2\NN}\sum_{j\in\ZZ}(-1)^j(\ch^jD\ch^h\SS_w^\sha,(i^w_h)^{(j)})
=u^{-l(w)-\nu}\sum_{h\in2\NN}(\ch^h\SS_w^\sha,i^w_h)\in\bK(\cc)$$
that is 
$$\DD(\sum_{h\in2\NN}(\ch^h\SS_w^\sha,i^w_h)=u^{-l(w)-\nu}\sum_{h\in2\NN}(\ch^h\SS_w^\sha,i^w_h)\in\bK(\cc).$$
Hence, setting
$$\fA_w=\sum_{h\in2\NN}(\ch^h\SS_w^\sha,i^w_h)\in\bK(\cc)$$
we have
$$\DD(\fA_w)=u^{-l(w)-\nu}\fA_w$$
that is
$$\ov{\fA_w}=u^{-l(w)}\fA_w$$
where we identify $\bK(\cc)=M$ as in 2.9. Under this identification the element $\fA_w\in M$ becomes 
$$\fA_w=\sum_{y\in\II;y\le w}P^\s_{y,w}a_y\in M.$$
From the definition of $\SS_w^\sha$ we have $\deg P^\s_{y,w}\le(l(w)-l(y)-1)/2$ (if $y\in\II,y<w$) and 
$P^\s_{w,w}=1$. We see that the existence part of 0.3(a) is verified by the element $A_w=v^{-l(w)}\fA_w\in\uM$ 
(recall that $v^2=u$).

\subhead 3.2\endsubhead
We prove the uniqueness part of 0.3(a). Assume that we have an element
$A'_w=v^{-l(w)}\sum_{y\in\II;y\le w}P'{}^\s_{y,w}a_y\in\uM$ 
($P'{}^\s_{y,w}\in\ZZ[u]$) such that $\ov{A'_w}=A'_w$, $P'{}^\s_{w,w}=1$ and for any $y\in\II$, $y<w$, we have 
$\deg P'{}^\s_{y,w}\le(l(w)-l(y)-1)/2$. We must show that $\cp_{z,w}=0$ where  
$\cp_{z,w}=P'{}^\s_{z,w}-P^\s_{z,w}$ for all $z\in\II,z\le w$.

We already know that $\cp_{w,w}=0$. We can assume that $z<w$ and that $\cp_{y,w}=0$ for any $y\in\II$ such that 
$z<y\le w$. With the notation in 0.2(b) we have
$$v^{l(w)}\sum_{y\in\II;y\le w}\ov{\cp_{y,w}}\sum_{y'\in\II;y'\le y}r_{y',y}a_{y'}=
v^{-l(w)}\sum_{y\in\II;y\le w}\cp_{y,w}a_y$$
hence 
$$v^{l(w)}\sum_{y\in\II;z\le y\le w}\ov{\cp_{y,w}}r_{z,y}=v^{-l(w)}\cp_{z,w}.$$
Using our inductive assumption this becomes
$$v^{l(w)}\ov{\cp_{z,w}}r_{z,z}=v^{-l(w)}\cp_{z,w}.$$
Using $r_{z,z}=v^{-2l(z}$ this becomes
$$v^{l(w)-l(z)}\ov{\cp_{z,w}}=v^{-l(w)+l(z)}\cp_{z,w}.$$
Here the right hand side is in $v\i\ZZ[v\i]$ and the left hand side is in $v\ZZ[v]$ (we use that $z<w$). Hence 
both sides are zero. Thus $\cp_{z,w}=0$. This completes the proof of 0.3(a). Now 0.3(b) is immediate. 

In the course of this proof we have also verified the following result.

\proclaim{Proposition 3.3} For any $y\in\II$, $y\le w$, the polynomial $P^\s_{y,w}$ defined in 0.3 satisfies
$$P^\s_{y,w}=\sum_{h\in2\NN}\tr((k^w_h)_y:V^{\ch^h\SS_w^\sha}_y@>>>V^{\ch^h\SS_w^\sha}_y)u^{h/2}\in\ZZ[u].$$
\endproclaim
Note that, by \cite{\KLL}, the polynomial $P_{y,w}$  of \cite{\KL} satisfies (for $y,w$ as in the proposition):
$$P_{y,w}=\sum_{h\in2\NN}\dim(V^{\ch^h\SS_w^\sha}_y)u^{h/2}\in\ZZ[u].$$

\subhead 3.4\endsubhead
Let $s\in S$ and let $w\in\II$ be such that $sw<w$ or equivalently $ws<w$. Let $y\in\II$ be such that $y\le w$.
Then we have also $sy\le w$; moreover, if $sy\ne ys$ we have $sys\le w$. We show:

(a) {\it If $sy=ys$ then $P^\s_{y,w}=P^\s_{sy,w}$.}

(b) {\it If $sy\ne ys$ then $P^\s_{y,w}=P^\s_{sys,w}$.}
\nl
Let $\cp$ be the variety of parabolic subgroups $P$ of $G$ such that for any Borel subgroups $B,B'$ in $P$
we have $(B,B')\in\bco_s$ and for some Borel subgroups $B,B'$ in $P$ we have $(B,B')\in\co_s$. Let $\r:\cb@>>>\cp$
be the morphism $B\m P$ where $P\in\cp$ contains $B$. Let $\ti\r=\r\T\r:\cb\T\cb@>>>\cp\T\cp$. This map commutes 
with the diagonal actions of $G$ and $\ti\r\s=\s'\ti\r$ where $\s':\cp\T\cp@>>>\cp\T\cp$ is $(P,P')\m(P',P)$.
We have $\bco_w=\ti\r\i(X)$ where $X=\ti\r(\bco_w)$, a closed subvariety of $\cp\T\cp$. Let $X_0=\ti\r(\co_w)$, 
the unique open $G$-orbit in $X$. Let $K$ be the intersection cohomology complex of $X$ with coefficients in
$\ql$ (on $X_0$) extended by $0$ on $\cp\T\cp-X$. We have $\SS_w^\sha=\ti\r^*K$.
Let $Z=\ti\r(\co_y)$; we have $\s'(Z)=Z$. For any $h\in2\NN$ there is a canonical vector space $V'$ with 
a linear involutive $\s'$-action which is canonical isomorphic to the stalk of $\ch^hK$ at any $z\in Z$.
From the definitions we have canonically 
$V^{\ch^h\SS_w^\sha}_y=V'$ and $V^{\ch^h\SS_w^\sha}_{sy}=V'$ if $sy=ys$,
(resp. $V^{\ch^h\SS_w^\sha}_{sys}=V'$ if $sy\ne ys$); moreover under these identifications the operators
$(k^w_h)_y$ and $(k^w_h)_{sy}$ if $sy=ys$ (resp. $(k^w_h)_{sys}$ if $sy\ne ys$) correspond to the
$\s'$-action on $V'$. It follows that 
$$\align&\tr((k^w_h)_y:V^{\ch^h\SS_w^\sha}y@>>>V^{\ch^h\SS_w^\sha}_y)=\tr(\s':V'@>>>V')\\&=
\tr((k^w_h)_{sy}:V^{\ch^h\SS_w^\sha}_{sy}@>>>V^{\ch^h\SS_w^\sha}_{sy})\endalign$$
if $sy=ys$ and
$$\align&\tr((k^w_h)_y:V^{\ch^h\SS_w^\sha}y@>>>V^{\ch^h\SS_w^\sha}_y)=\tr(\s':V'@>>>V')\\&=
\tr((k^w_h)_{sys}:V^{\ch^h\SS_w^\sha}_{sys}@>>>V^{\ch^h\SS_w^\sha}_{sys})\endalign$$
if $sy\ne ys$. This implies (a) in view of 3.3.

\head 4. The action of $u\i(T_s+1)$ in the basis $(A_w)$ \endhead
\subhead 4.1\endsubhead
Let $y,w\in\II$. We set $\d_{y,w}=1$ if $y=w$, $\d_{y,w}=0$ if $y\ne w$, $\d'_{y,w}=1-\d_{y,w}$.
When $y\le w$ we set $\p_{y,w}=v^{-l(w)+l(y)}P^\s_{y,w}$ so that $\p_{y,w}\in v\i\ZZ[v\i]$ if 
$y<w$ and $\p_{w,w}=1$; when $y\not\le w$ we set $\p_{y,w}=0$. In any case we have

(a) $\p_{y,w}=\d_{y,w}+\mu'(y,w)v\i+\mu''(y,w)v^{-2}\mod v^{-3}\ZZ[v\i]$
\nl
where $\mu'(y,w)\in\ZZ,\mu''(y,w)\in\ZZ$. Note that

(b) $\mu'(y,w)\ne0\implies y<w, l(y)\ne l(w)\mod2$,

(c) $\mu''(y,w)\ne0 \implies y<w, l(y)=l(w)\mod2$.
\nl
For $w\in\II$ we set $a'_w=v^{-l(w)}a_w$ so that $A_w=\sum_{y\in\II}\p_{y,w}a'_y$.

\subhead 4.2\endsubhead
In this section we fix $s\in S$. We set  $c_s=v^{-2}(T_s+1)\in\ufH'$. The formulas in 0.2(a) (with $w\in\II$) can
be rewritten as follows.

$c_sa'_w=(v+v\i)a'_{sw}+(1+v^{-2})a'_w$ if $sw=ws>w$;

$c_sa'_w=(v-v\i)a'_{sw}+(v^2-1)a'_w$ if $sw=ws<w$;

$c_sa'_w=a'_{sws}+v^{-2}a'_w$ if $sw\ne ws>w$;

$c_sa'_w=a'_{sws}+v^2a'_w$ if $sw\ne ws<w$.

\subhead 4.3\endsubhead
For any $y,w\in\II$ such that $sy<y<sw>w$ we define $\cm^s_{y,w}\in\uca$ by:
$$\cm^s_{y,w}=\mu''_{y,w}-\sum_{x\in\II;y<x<w,sx<x}\mu'_{y,x}\mu'_{x,w}
-\d_{w,sws}\mu'_{y,sw}+\mu'_{sy,w}\d_{sy,ys}$$
if $l(y)=l(w)\mod2$,
$$\cm^s_{y,w}=\mu'_{y,w}(v+v\i)$$ 
if $l(w)\ne l(y)\mod2$.

\proclaim{Theorem 4.4} Recall that $s\in S$. Let $w\in\II$. 

(a) If $sw=ws>w$ then $c_sA_w=(v+v\i)A_{sw}+\sum_{z\in\II;sz<z<sw}\cm^s_{z,w}A_z$.

(b) If $sw\ne ws>w$ then $c_sA_w=A_{sws}+\sum_{z\in\II;sz<z<sw}\cm^s_{z,w}A_z$.

(c) If $ws<w$ then $c_sA_w=(v^2+v^{-2})A_w$.
\endproclaim
We prove (c). For $z\in\II,sz<z$ we set $\ti a'_z=a'_z+v\i a'_{sz}$ (if $sw=ws$) and 
$\ti a'_z=a'_z+v^{-2}a'_{szs}$ (if $sz\ne zs$). By 3.4 we have 
$$A_z=\ti a'_z+\text{$\uca$-linear combination of elements $\ti a'_y$ with $y\in\II$, $sy<y<z$}.\tag d$$
It follows that the elements $A_z (z\in\II,sz<z)$ span the same $\uca$-submodule $\uM'$ of $\uM$ as the elements
$\ti a'_z (z\in\II,sz<z)$. To show (c) it is enough to show that $c_s-(v^2+v^{-2})$ acts as $0$ on $\uM'$; hence
it is enough to show that $(c_s-(v^2+v^{-2}))\ti a'_z=0$ for any $z\in\II,sz<z$. But this follows from 4.2. This 
proves (c).

In the rest of the proof we assume that $sw>w$. Note that 
$$c_sA_w=\sum_{y\in\WW;y\le w}\p_{y,w}c_sa'_y$$
hence using 4.2, $c_sA_w$ is an $\uca$-linear combination
of elements of the form $\ti a'_z$  with $sz<z$ and with $z\le sw$ (if $sw=ws$) or with $z\le sws$ (if 
$sw\ne ws$). Using (d) it follows that $c_sA_w$ is an $\uca$-linear combination of elements of the form $A_z$ 
with $sz<z$ and with $z\le sw$ (if $sw=ws$) or with $z\le sws$ (if $sw\ne ws$). For such $z$ we have either 
$z=sw$ or $z=sws$ or $z<sw$.
(To see this we can assume that $sw\ne ws$, $z<sws$. Then we must have $z\le sw$ or $z\le ws$; 
if $z\le ws$ then taking inverses we obtain $z\le sw$. Thus $z\le sw$ in any case. We cannot have $z=sw$ since
$sw\ne ws$. hence $z<sw$.) We see that $c_sA_w=\sum_{x\in\II;sx<x}m_xA_x$ where $m_x\in\uca$; moreover we have 
$m_x=0$ unless $x=sw$ (if $sw=ws$), $x=sws$ (if $sw\ne ws$) or $x<sw$.

Since $\ov{c_sA_w}=c_sA_w$ we have $\sum_{x\in\II;sx<x}\bar{m_x}A_x=\sum_{x\in\II;sx<x}m_xA_x$ hence 
$\ov{m_x}=m_x$ for all $x\in\II$ such that $sx<x$. We have
$$\align&c_sA_w=\sum_{y\in\II;sy=ys>y}\p_{y,w}c_sa'_y+\sum_{y\in\II;sy=ys<y}\p_{y,w}c_sa'_y\\&
+\sum_{y\in\II;sy\ne ys>y}\p_{y,w}c_sa'_y+\sum_{y\in\II;sy\ne ys<y}\p_{y,w}c_sa'_y\endalign$$
hence, using 4.2:
$$\align&c_sA_w=\sum_{y\in\II;sy=ys>y}\p_{y,w}((v+v\i)a'_{sy}+(1+v^{-2})a'_y)\\&+
\sum_{y\in\II;sy=ys<y}\p_{y,w}((v-v\i)a'_{sy}+(v^2-1)a'_y)\\&+
\sum_{y\in\II;sy\ne ys>y}\p_{y,w}(a'_{sys}+v^{-2}a'_y)+\sum_{y\in\II;sy\ne ys<y}\p_{y,w}(a'_{sys}+v^2a'_y),
\endalign$$
$$\align&c_sA_w=\sum_{y\in\II;sy=ys<y}\p_{sy,w}(v+v\i)a'_y+\sum_{y\in\II;sy=ys>y}\p_{y,w}(1+v^{-2})a'_y\\&
+\sum_{y\in\II;sy=ys>y}\p_{sy,w}(v-v\i)a'_y+\sum_{y\in\II;sy=ys<y}\p_{y,w}(v^2-1)a'_y\\&
+\sum_{y\in\II;sy\ne ys<y}\p_{sys,w}a'_y+\sum_{y\in\II;sy\ne ys>y}\p_{y,w}v^{-2}a'_y\\&
+\sum_{y\in\II;sy\ne ys>y}\p_{sys,w}a'_y+\sum_{y\in\II;sy\ne ys<y}\p_{y,w}v^2a'_y.\endalign$$
Thus,
$$\align&c_sA_w=\sum_{y\in\II;sy=ys<y}(\p_{sy,w}(v+v\i)+\p_{y,w}(v^2-1))a'_y\\&
+\sum_{y\in\II;sy=ys>y}(\p_{y,w}(1+v^{-2})+\p_{sy,w}(v-v\i))a'_y\\&
+\sum_{y\in\II;sy\ne ys<y}(\p_{sys,w}+\p_{y,w}v^2)a'_y+\sum_{y\in\II;sy\ne ys>y}(\p_{y,w}v^{-2}+\p_{sys,w})a'_y.
\endalign$$
We have
$$\sum_{x\in\II;sx<x}m_xA_x=\sum_{y\in\II}\sum_{x\in\II;sx<x}m_x\p_{y,x}a'_y.$$
It follows that for $y\in\II$ such that $ys<y$ we have
$$\sum_{x\in\II;sx<x}m_x\p_{y,x}=\p_{sy,w}(v+v\i)+\p_{y,w}(v^2-1)\text{ if }sy=ys;$$
$$\sum_{x\in\II;sx<x}m_x\p_{y,x}=\p_{sys,w}+\p_{y,w}v^2\text{ if }sy\ne ys.$$
For any $f\in\uca$ we define $f^+\in\ZZ[v]$ by $f-f^+\in v\i\ZZ[v\i]$. It follows that 
$$\sum_{x\in\II;sx<x}(m_x\p_{y,x})^+=(\p_{sy,w}(v+v\i))^++(\p_{y,w}(v^2-1))^+\text{ if }sy=ys<y;$$
$$\sum_{x\in\II;sx<x}(m_x\p_{y,x})^+=\p_{sys,w}^++(\p_{y,w}v^2)^+\text{ if }sy\ne ys<y.$$
Using 4.1(a) and that $y\ne w$ if $sy<y$ we deduce
$$\align&m_y^++\sum_{x\in\II;sx<x,y<x}(m_x\p_{y,x})^+=\mu'_{sy,w}+\d_{sy,w}v+\mu''_{y,w}+\mu'_{y,w}v
\text{ if }sy=ys<y;\\&
m_y^++\sum_{x\in\II;sx<x,y<x}(m_x\p_{y,x})^+=\d_{sys,w}+\mu''_{y,w}+\mu'_{y,w}v\text{ if }sy\ne ys<y.\tag e
\endalign$$
In particular we have 
$$m_y^++\sum_{x\in\II;sx<x,y<x}(m_x\p_{y,x})^+\in\ZZ+\ZZ v$$ 
for any $y\in\II$ such that $sy<y$. This shows by descending induction on $l(y)$ that $m_y^+\in\ZZ+\ZZ v$ for any
$y\in\II$ such that $sy<y$. (Indeed, if we know that for $x\in\II$ such that $sx<x,y<x$ we have 
$m_x^+\in\ZZ+\ZZ v$, then $(m_x\p_{y,x})^+\in\ZZ$.) Setting $m_y^+=m^0_y+m'_yv$ (with $m^0_y\in\ZZ$, 
$m'_y\in\ZZ$) for any $y\in\II$ such that $sy<y$, we can rewrite (e) as follows: 
$$\align&m_y^0+m'_yv+\sum_{x\in\II;sx<x,y<x}m'_x\mu'_{y,x}\\&=
\d_{sy,ys}\mu'_{sy,w}+\d_{sy,ys}\d_{sy,w}v+\d'_{sy,ys}\d_{sys,w}+\mu''_{y,w}+\mu'_{y,w}v.\endalign$$
In particular,
$$m'_y=\mu'_{y,w}+\d_{sy,ys}\d_{sy,w}$$
for any $y\in\II$ such that $sy<y$. It follows that 
$$\align&m_y^0+\sum_{x\in\II;sx<x,y<x}(\mu'_{x,w}+\d_{sx,xs}\d_{sx,w})\mu'_{y,x}\\&=
\d_{sy,ys}\mu'_{sy,w}+\d'_{sy,ys}\d_{sys,w}+\mu''_{y,w}.\endalign$$
Equivalently we have
$$\align&m_y^0=-\sum_{x\in\II;sx<x,y<x}\mu'_{y,x}\mu'_{x,w}-\d_{w,sws}\mu'_{y,sw}\\&
+\d_{sy,ys}\mu'_{sy,w}+\d'_{sy,ys}\d_{sys,w}+\mu''_{y,w}.\endalign$$
(We have used that $\sum_{x\in\II;sx<x,y<x}\d_{sx,xs}\d_{sx,w}\mu'_{y,x}=\d_{w,sws}\mu'_{y,sw}$.)
Since $\ov{m_y}=m_y$ we must have $m_y=m_y^0+m'_y(v+v\i)$ for $y\in\II,sy<y$.
For $y\in\II$ such that $sy<y,l(w)\ne l(y)\mod2$ we deduce using 4.1(c) that
$$m_y=\mu'_{y,w}(v+v\i)=\cm^s_{y,w}\text{ if $y<sw$},\qua m_y=v+v\i$$
if $y=sw$, hence $ws=sw$. For $y\in\II$ such that $sy<y,l(w)=l(y)\mod2$ we deduce using 4.1(b),(c) that
$$m_y=\mu''_{y,w}-\sum_{x\in\II;sx<x,y<x}\mu'_{y,x}\mu'_{x,w}-\d_{w,sws}\mu'_{y,sw}
+\d_{sy,ys}\mu'_{sy,w}+\d'_{sy,ys}\d_{sys,w}.$$
If $ws\ne sw$ and $y=sws$ we deduce that $m_y=1$. If $y\in\II,sy<y,y<sw,l(y)=l(w)\mod2$ then
$$m_y=\mu''_{y,w}-\sum_{x\in\II;sx<x,y<x}\mu'_{y,x}\mu'_{x,w}-\d_{w,sws}\mu'_{y,sw}
+\d_{sy,ys}\mu'_{sy,w}=\cm^s_{y,w}.$$
We see that $m_y=\cm^s_{y,w}$ for any $y\in\II$ such that $sy<y<sw$. We also see that $m_{sw}=v+v\i$ if $sw=ws$ 
and $m_{sws}=1$ if $sw\ne ws$. This completes the proof of the theorem.

\subhead 4.5\endsubhead
We now present an algorithm for compute the polynomials $P^\s_{y,w}$ for $y\le w$ in $\II$. It will be convenient
to state this in terms of the elements $\p_{y,w}=v^{-l(w)+l(y)}P^\s_{y,w}\in\ZZ[v\i]$ (see 4.1). Recall that
$\p_{y,w}$ is defined to be $0$ if $y,w\in\II,y\not\le w$.

We can restate Theorem 4.4(a),(b) as follows. (Recall that $s\in S$, $w\in\II$, $sw>w$.)
$$\align&\sum_{y\in\II;sy=ys<y}(\p_{sy,w}(v+v\i)+\p_{y,w}(v^2-1))a'_y\\&
+\sum_{y\in\II;sy=ys>y}(\p_{y,w}(1+v^{-2})+\p_{sy,w}(v-v\i))a'_y\\&
+\sum_{y\in\II;sy\ne ys<y}(\p_{sys,w}+\p_{y,w}v^2)a'_y+\sum_{y\in\II;sy\ne ys>y}(\p_{y,w}v^{-2}+\p_{sys,w})a'_y\\&
=\sum_{y\in\II}\sum_{x\in\II;sx<x<sw}\cm^s_{x,w}\p_{y,x}a'_y\\&
+\sum_{y\in\II}(v+v\i)\p_{y,sw}\d_{sw,ws}a'_y+\sum_{y\in\II}\p_{y,sws}(1-\d_{sw,ws})a'_y.\endalign$$
It follows that for any $y\in\II$, the expression
$$(v+v\i)\p_{y,sw}\d_{sw,ws}+\p_{y,sws}\d'_{sw,ws}$$
is equal to
$$-\sum_{x\in\II;sx<x<sw}\cm^s_{x,w}\p_{y,x}+\p_{sy,w}(v+v\i)+\p_{y,w}(v^2-1)\text{ if }sy=ys<y,$$
$$-\sum_{x\in\II;sx<x<sw}\cm^s_{x,w}\p_{y,x}+\p_{y,w}(1+v^{-2})+\p_{sy,w}(v-v\i)\text{ if }sy=ys>y,$$
$$-\sum_{x\in\II;sx<x<sw}\cm^s_{x,w}\p_{y,x}+\p_{sys,w}+\p_{y,w}v^2\text{ if }sy\ne ys<y,$$
$$-\sum_{x\in\II;sx<x<sw}\cm^s_{x,w}\p_{y,x}+\p_{y,w}v^{-2}+\p_{sys,w}\text{ if }sy\ne ys>y.$$
We want show that the formulas above determine uniquely the quantities $\p_{y,sw}$ (resp. $\p_{y,sws}$)
assuming that $y\in\II$ and that $sw=ws$ (resp. $sw\ne ws$) and assuming that the
quantities $\p_{y',w'}$ are known for any $w'\in\II$ such that $l(w')<l(ws)$ (resp. $l(w')<l(sws)$)
and any $y'\in\II$. Then the quantities $\cm^s_{x,w}$ in these formulas are also known except for
a part of them given by $\d_{ws,sw}\mu'_{x,sw}$ which is not known.
If in the formulas above we replace the terms that are assumed to be known by a symbol $\sp$ we obtain
$$(v+v\i)\p_{y,sw}=\sum_{x\in\II;sx<x<sw}\mu'_{x,sw}\p_{y,x}+\sp\text{ if $ws=sw$}.\tag a$$
$$\p_{y,sws}=\sp\text{ if }sw\ne ws.$$
We can now assume that $sw=ws$. 
In this case we determine the quantities $\p_{y,sw}$ by descending induction on $l(y)$. (We can assume that
$y<sw$ since $\p_{y,sw}=1$ for $y=sw$.) Thus we can
assume that $\p_{y,sw}$ (hence also $\mu'_{y,sw}$) is known when $y$ is replaced by $x\in\II$ such that $y<x<sw$.
Since in the sum over $x$ in (a) we can restrict to those $x$ such that $y\le x$ we see that (a) becomes
$$(v+v\i)\p_{y,sw}-\mu'(y,sw)=\sum_{x\in\II;sx<x<sw,y<x}\mu'_{x,sw}\p_{y,x}+\sp$$
that is
$$(v+v\i)\p_{y,sw}-\mu'(y,sw)=\sp.$$
Let us write $\p_{y,sw}=\sum_{n\ge1}c_nv^{-n}$ where $c_n\in\ZZ$ are zero for all but finitely many $n$; note
that $c_1=\mu'(y,sw)$. (Recall that $y<sw$.) It follows that
$$\sum_{n\ge1}c_nv^{-n+1}+\sum_{n\ge1}c_nv^{-n-1}-c_1=\sp$$
so that $c_2=\sp$, $c_1+c_3=\sp$, $c_2+c_4=\sp$, $c_3+c_5=\sp$, $c_4+c_6=\sp$, $\do$.
It follows that $c_{2k}=\sp$ for $k=1,2,\do$ and, since $c_{2t+1}=0$ for large 
$t$ we have also $c_{2k-1}=\sp$ for $k=1,2,\do$. Thus $c_i=\sp$ for $i=1,2,\do$ so that $\p_{y,sw}=\sp$.

The procedure above gives an algorithm to compute $\p_{y,z}$ for any $y,z\in\II$ such that $y\le z$. Indeed, if 
$z=1$ then $y=1$ and $\p_{y,z}=1$. If $z\ne1$ then we can find $s\in S$ such that $sz<z$. Setting $w=zs$ (if 
$zs=sz$) or $w=szs$ (if $zs\ne sz$) we see that $\p_{y,z}$ is determined by the inductive procedure above.

\head 5. Relation with two-sided cells\endhead
\subhead 5.1\endsubhead
For any $w\in W$ let $\dot{c}_w=v^{-l(w)}\sum_{y\in W;y\le w}P_{y,w}(v^2)T_y\in\ufH$ (compare \cite{\KL}); 
similarly let $c_w=u^{-l(w)}\sum_{y\in W;y\le w}P_{y,w}(u^2)T_y\in\ufH'$. The elements $\dot{c}_w (w\in W)$ form
an $\uca$-basis of $\ufH$; the elements $c_w (w\in W)$ form an $\ca$-basis of $\ufH'$. 
For $z,w$ in $W$ we write (in the algebra $\ufH$):

$\dot{c}_z\dot{c}_w\dot{c}_{z\i}=\sum_{w'\in W}h_{z,w,w'}\dot{c}_{w'}$ 
\nl
where $h_{z,w,w'}\in\NN[v,v\i]$. For $z\in W,w\in\II$ we write (using the $\ufH'$-module structure on $\uM$):

$c_zA_w=\sum_{w'\in\II}f_{z,w,w'}A_{w'}$
\nl
where $f_{z,w,w'}\in\uca$ are related to $h_{z,w,w'}$ as follows.
For $z\in W$, $w,w'\in\II$ we write $h_{z,w,w'}=\sum_{n\in\ZZ}b_nv^n$, $f_{z,w,w'}=\sum_{n\in\ZZ}b'_nv^n$
where $b_n\in\NN,b'_n\in\ZZ$ are zero for all but finitely many $n$. One can show that for each $n$ we have 
$b_n=b_n^++b_n^-$, 
$b'_n=b_n^+-b_n^-$ for some $b_n^+\in\NN,b_n^-\in\NN$. (This is similar to the relation between $P_{y,w}$ and 
$P^\s_{y,w}$ for $y,w\in\II$. It is also analogous to the phenomenon described in \cite{\UN, 16.3(a),(b)}.) In 
particular,

(a) {\it If for some $n$ we have $b_n=0$, then $b'_n=0$. Hence if $f_{z,w,w'}\ne0$, then $h_{z,w,w'}\ne0$.}

(b) {\it If for some $n$ we have $b_n=1$, then $b'_n=\pm1$.}

\subhead 5.2\endsubhead
Let $c$ be a two-sided cell of $W$, see \cite{\KL}.
For $y,w\in W$ we shall write $y\preceq w$ instead of $y\le_{LR}w$ ($\le_{LR}$ is the preorder defined in
\cite{\KL}). For $y\in W$ we write $y\preceq c$ instead of: $y\preceq w$ for some $w\in c$.
For $y\in W$ we write $y\prec c$ instead of: $y\preceq c$ and $y\n c$.

Let $\ufH'{}^{\preceq c}$ be the $\uca$-submodule of $\ch'$ spanned by the elements $c_{w'}$ 
where $w'\in W$ is such that $w'\preceq c$. Let $\ufH'{}^{\prec c}$ be the $\uca$-submodule of $\ufH'$ spanned 
by the elements $c_{w'}$ where $w'\in W$ is such that $w'\prec c$.
Note that $\ufH'{}^{\preceq c},\ufH'{}^{\prec c}$ are two-sided ideals of $\ufH'$. Hence 
$\ufH'{}^{\preceq c}/\ufH'{}^{\prec c}$ is naturally an $\ufH'$-bimodule, in particular a left $\ufH'$-module.

Let $\uM^{\preceq c}$ be the $\uca$-submodule of $\uM$ spanned by
the elements $A_{w'}$ where $w'\in\II$, $w'\preceq c$. Let $\uM^{\prec c}$ be the
$\uca$-submodule of $\uM$ spanned by the elements $A_{w'}$ where $w'\in\II$, $w'\prec c$.

From 5.1(a) we see that for $z\in W$, $w\in\II$, the element $c_zA_w$ is an $\uca$-linear combination of elements
$A_{w'}$ with $w'\in\II,w'\preceq w$. In particular, $\uM^{\preceq c},\uM^{\prec c}$ are $\ufH'$-submodules of 
$uM$. Hence $\uM^{\preceq c}/\uM^{\prec c}$ is naturally an $\ufH'$-module.

Let $\a\in\NN$ be the value of the $\aa$-function \cite{\UN, \S13} on $c$. 
Let $\et\in\II\cap c$. We define an $\ufH'$-linear map
$\tau_\et:\ufH'{}^{\preceq c}@>>>\uM^{\preceq c}$ by $\x\m v^{2\a}\x A_\et$ (we use the $\ufH'$-module structure 
on $\uM$). From 5.1(a) we see that if $\x\in\ufH'{}^{\prec c}$ then $\x A_\et\in\uM^{\prec c}$. Thus $\tau_\et$ 
restricts to an $\ufH'$-linear map $\ufH'{}^{\prec c}@>>>\uM^{\prec c}$ hence it induces
an $\ufH'$-linear map $\ufH'{}^{\preceq c}/\ufH'{}^{\prec c}@>>>\uM^{\preceq c}/\uM^{\prec c}$ denoted again by 
$\tau_\et$.
Let $d$ be the unique distinguished involution in the same left cell as $\et$ (hence in the same right cell as 
$\et$). From the known properties of distinguished involutions (see \cite{\UN, \S15}) we see that the 
following holds in $\ufH$.
$$v^{2\a}\dot{c}_d\dot{c}_\et\dot{c}_d=\dot{c}_\et+\x+\x'$$
where $\x\in\sum_{x\in c}v\ZZ[v]\dot{c}_x$ and $\x'\in\ufH^{\prec c}$.
Using now the results in 5.1 we deduce that
$$v^{2\a}c_dA_\et=A_\et+m_\et+m'$$
where $m_\et\in\sum_{x\in c}v\ZZ[v]A_x$ and $m'\in\uM^{\prec c}$. Thus we have
$$\tau_\et(c_d)=\pm A_\et+m_\et$$
where $m_\et$ is as above. Now let $\tau:\op_{\et\in c}\ufH'{}^{\preceq c}/\ufH'{}^{\prec c}@>>>
\uM^{\preceq c}/\uM^{\prec c}$ 
be the $\ufH'$-linear map whose restriction to the $\et$-summand is $\tau_\et$. The image of this map contains the
elements $\pm A_\et+m_\et$ ($\et\in c$) which clearly form a basis of 
$\ZZ[[v]]\ot_{\uca}(\uM^{\preceq c}/\uM^{\prec c})$: the $c\T c$ matrix whose $\et,\et'$-entry (in $\ZZ[v]$) 
is the $A_{\et'}$-coordinate of $\pm A_\et+m_\et$ has determinant $\pm1$ plus an element in $v\ZZ[v]$, hence is
invertible in $\ZZ[[v]]$. Thus after extension of scalars to $\ZZ[v]$, $\tau$ is surjective. Hence after 
extension of scalars to the quotient field of $\ZZ[[v]]$, or to the quotient field $\QQ(v)$ of $\uca$, $\tau$
is surjective. We see that 

(a) {\it the $\QQ(v)\ot_{\uca}\ufH'$-module $\QQ(v)\ot_{\uca}(\uM^{\preceq c}/\uM^{\prec c})$ is a direct
sum of irreducible $\QQ(v)\ot_{\uca}\ufH'$-modules which appear in the $\QQ(v)\ot_{\uca}\ufH'$-module 
$\QQ(v)\ot_{\uca}(\ufH'{}^{\preceq c}/\ufH'{}^{\prec c})$ carried by the two-sided cell $c$.}

\subhead 5.3\endsubhead
In this subsection we assume that $W$ is a Weyl group of type $A_{n-1}$, $n\ge2$. In this case for any two sided 
cell $c$, the $\QQ(v)\ot_{\uca}\ufH'$-module $\QQ(v)\ot_{\uca}(\ufH'{}^{\preceq c}/\ufH'{}^{\prec c})$ is known 
to be a direct sum of copies of a single simple module $E_c$. Hence from 5.2(a) we deduce that the 
$\QQ(v)\ot_{\uca}\ufH'$-module $\QQ(v)\ot_{\uca}(\uM^{\preceq c}/\uM^{\prec c})$ is a direct sum of say $n_c$
copies of $E_c$. The dimension over $\QQ(v)$ of this module is the number of involutions contained in $c$ which
is known to be equal to $\dim E_c$. It follows that $n_c=1$. It follows that the $\QQ(v)\ot_{\uca}\ufH'$-module 
$\QQ(v)\ot_{\uca}\uM$ is isomorphic to $\op_cE_c$ that is, is a ``model representation" for the algebra
$\QQ(v)\ot_{\uca}\ufH'$. (Each simple module of the algebra appears exactly once in it.)

\head 6. The $W$-module $M_1$\endhead
\subhead 6.1\endsubhead
Let $\un{\II}$ be the set of conjugacy classes of $W$ contained in $\II$. For any $w\in\II$ let $h(w)$ be the 
dimension of the fixed point set $R^{-w}$ of $-w$ in the reflection representation $R$ of $W$ (over $\QQ$). For 
any $C\in\un{\II}$ we set $h(C)=h(w)$ where $w$ is any element of $C$. We have the following result.

\proclaim{Lemma 6.2} Let $s\in S$, $w\in\II$ be such that $sw=ws>w$. Then $h(sw)>h(w)$.
\endproclaim
From our assumptions it follows that the line $R^{-s}$ is contained in $R^w$, the fixed point set of $w$, and 
$R^{-w}$ is contained in $R^s$, the fixed point set of $s$. It follows that $R^{-s}$ and $R^{-w}$ are contained
in $R^{-sw}$ hence $R^{-s}\op R^{-w}$ is contained in $R^{-sw}$. Thus $h(sw)\ge h(w)+1$. The lemma is proved.

\subhead 6.3\endsubhead
In the setup of 0.4 for any $i\in\NN$ let $M_1^{\ge i}$ be the subspace of $M_1$ spanned by $\{a_w;h(w)\ge i\}$.
From the formulas for the $W$-action in 0.4 and from 6.2 we see that for any $i\in\NN$, $M_1^{\ge i}$ is a 
$W$-submodule of $M_1$. This induces a $W$-module structure on $M_1^i:=M_1^{\ge i}/M_1^{\ge i+1}$.
Let $gr(M_1)=\op_{i\in\NN}(M_1^{\ge i}/M_1^{\ge i+1})$. This is naturally a $W$-module. By complete reducibility
we have $M_1\cong gr(M_1)$ as $W$-modules. For any $w\in\II$ we denote by $\ti a_w$ the image of
$a_w$ under the projection $M_1^{\ge h(w)}@>>>M_1^{\ge h(w)}/M_1^{\ge h(w)+1}$. The elements $\ti a_w (w\in\II)$
form a $\QQ$-basis of $gr(M_1)$ in which the $W$-action is as follows. (Here $s\in S$.)

$s(\ti a_w)=-\ti a_{sws}$ if $w\in\II,sw=ws<w$;

$s(\ti a_w)=\ti a_{sws}$ for all other $w\in\II$.
\nl
It folows that for any $x\in W,w\in\II$ we have $x(\ti a_w)=\e_{x,w}\ti a_{xwx\i}$ where $\e_{x,w}\in\{1,-1\}$
satisfies $\e_{xy,w}=\e_{x,ywy\i}\e_{y,w}$ for all $x,y\in W,w\in\II$. In particular, if $x,y$ are in $Z(w)$, the
centralizer of $w\in\II$, in $W$ then $\e_{xy,w}=\e_{x,w}\e_{y,w}$.
Thus, $x\m\e_{x,w}$ is a homomorphism $\e_w:Z(w)@>>>\{1,-1\}$ and it is clear that for any $C\in\un{\II}$, the 
subspace of $gr(M_1)$ spanned by $\{\ti a_w;w\in C\}$ is a $W$-submodule of $gr(M_1)$ isomorphic to 
the representation of $W$ induced from the character $\e_w$ of $Z(w)$ where $w$ is any element of $C$. We see that

(a) $M_1\cong\op_w\Ind_{Z(w)}^W(\e_w)$
\nl
as a $W$-module; here $w$ runs over a set of representatives for the various $C\in\un{\II}$.
The multiplicities of the various irreducible $W$-modules in the $W$-module given by the right hand side of (a) 
were first described explicitly in \cite{\IRS} for $W$ of type $A_{n-1}$ and later by Kottwitz \cite{\KO} for any
irreducible $W$; for example, if $W$ is irreducible of classical type then an irreducible representation $E$ of 
$W$ appears in
the right hand side of (a) if and only if it is special in the sense of \cite{\OR, (4.1.4)} and then its 
multiplicity is the integer $f_E$ in \cite{\OR, (4.1.1)} (in type $A_{n-1}$ this follows also from 5.3). 

\subhead 6.4\endsubhead
Kottwitz's formula for the $W$-module $gr M_1$ in the right hand side of 6.3(a) is in terms of the 
nonabelian Fourier transform matrix of \cite{\OR, 4.14}. We shall reformulate that formula in terms of
unipotent representations, by using \cite{\OR, 4.23} which relates unipotent representations with 
the nonabelian Fourier transform matrix. One advantage of this reformulation is that 
unlike in \cite{\KO} we need not consider separately the ``exceptional" representations of $W$ of type
$E_7,E_8$. (In the remainder of this section we assume that $G$ is almost simple and that $q,\ph'$ are as in 1.1.)
 Let $Irr G^{\ph'}$ be a set of representatives for the isomorphism classes of
irreducible representations of $G^{\ph'}$ over $\bar\ql$. For any $\r\in Irr G^{\ph'}$ 
let $\e(\r)\in\{0,1,-1\}$ be the
Frobenius-Schur indicator of $\r$; thus $\e(\r)=1$ (resp. $\e(\r)=-1$) if $\r$ admits a nondegenerate 
$G^{\ph'}$-invariant symmetric (resp. antisymmetric) bilinear form $\r\T\r@>>>\bar\ql$ and $\e(\r)=0$ if
$\r$ is not isomorphic to its dual representation. It is known (Frobenius-Schur) that the dimension
of the virtual representation $\sum_{\r\in\Irr G^{\ph'}}\e(\r)\r$ is equal to the number of $g\in G^{\ph'}$ such
that $g^2=1$. Let us now consider the part 
$\Th:=\sum_{\r\in\cu}\e(\r)\r$ of the virtual representation above coming from the set
$\cu$ of unipotent representations of $G^{\ph'}$.

Let $\cx$ be the set of all triples 
$(\cf,y,r)$ where $\cf$ is a family \cite{\OR, 4.2} of irreducible representations of $W$ (with an associated 
finite group $\cg_\cf$, see \cite{\OR, Ch.4}), $y$ is an element of $\cg_\cf$ defined up to conjugacy and $r$ is 
an irreducible representation of the centralizer of $y$ in $\cg_\cf$ defined up to isomorphism. In 
\cite{\OR, 4.23}, $\cx$ is put in a bijection $(\cf,y,r)\lra\r_{\cf,y,r}$ with $\cu$. If $(\cf,y,r)\in\cx$ then 
we have $\e(\r_{\cf,y,r})\in\{0,1\}$; moreover, $\e(\r_{\cf,y,r})=1$ in exactly the following cases:

(i) $|\cf|\ne2$ and $y$ acts on $r$ by the scalar $\pm1$;

(ii) $|\cf|=2$ and $y=1$.
\nl
(This follows from results in \cite{\RA, Sec.1} where it is shown that in case (i) $\r_{\cf,y,r}$ is actually 
defined over $\QQ$. On the other hand, in case (ii), $\r_{\cf,y,r}$ is defined over $\QQ[\sqrt{q}]$.)
In particular, if $W$ is of classical type, we have 
$\e(\r_{\cf,y,r})=1$ for any $(\cf,y,r)\in\cx$.

For each $W$-module $E$ we define $R_E=|W|\i\sum_{w\in W}\tr(w,E)\sum_i(-1)^iH^i_c(X_w,\bar\ql)$ where $X_w$ is 
the variety associated to $G,\ph',w$ in \cite{\DL}. We view $R_E$ as an element of $\QQ[\cu]$, the vector space 
of $\QQ$-linear combinations of elements of $\cu$ with its symmetric inner product $(,)$ in which $\cu$ is an 
orthonormal basis.

Using the above description of $\e(\r_{\cf,y,r})$ and \cite{\OR, 4.23} we can rewrite Kottwitz's formula in the 
form

(a) {\it $\sum_{\r\in\cu;\e(\r)=1}(\r,R_E)=(R_{gr M_1},R_E)$ for any irreducible $W$-module $E$.}
\nl
Using the isomorphism $M_1\cong gr M_1$ we can rewrite (a) as follows.

(b) {\it We have $\Th=R_{M_1}+\xi$ where $\xi\in\QQ[\cu]$ is orthogonal to $R_{E}$ for any irreducible 
$W$-module $E$.}
\nl
Note that if $G$ is of type $\ne F_4,E_8$ then $(\xi,\xi)=0$ hence $\xi=0$; if $G$ is of type $F_4$ or $E_8$ then 
$(\xi,\xi)=1$.

\subhead 6.5\endsubhead
Let $c$ be a two-sided cell of $W$. Let $M^{\preceq c}_1$ be the subspace of $M_1$ spanned by the elements 
$A_{w'}$  where $w'\in\II$, $w'\preceq c$. Let $M^{\prec c}_1$ be the
subspace of $M_1$ spanned by the elements $A_{w'}$ where $w'\in\II$, $w'\prec c$. Then
$M^{\preceq c}_1,M^{\prec c}_1$ are $W$-submodules of $M_1$ and $M^{\preceq c}_1/M^{\prec c}_1$ is naturally a 
$W$-module. From 5.2(a) we see that this last $W$-module is a direct sum of irreducible representations $E$ 
of $W$ which appear in $c$ (which we write as $E\dsv c$). Since $M_1$ is isomorphic as a $W$-module to
$\op_{c'}M^{\preceq c'}_1/M^{\prec c'}_1$ ($c'$ runs over the two-sided cells of $W$) we see that
$M^{\preceq c}_1/M^{\prec c}_1\cong\op_{E;E\dsv c}E^{\op(E:M_1)}$ where $(E:M_1)$ is the multiplicity of $E$ in 
the
$W$-module $M_1$. Note that $\dim(M^{\preceq c}_1/M^{\prec c}_1)=|c\cap\II|$. It follows that 

(a) $|c\cap\II|=\op_{E;E\dsv c}(E:M_1)\dim(E)$.
\nl
This determines explicitly the number $|c\cap\II|$ since the multiplicities
$(E:M_1)=(E:gr M_1)$ are known from \cite{\KO}. In particular, if $W$ is of classical type we have

(b) $|c\cap\II|=f_E\dim(E)$
\nl
where $E$ is the special representation such that $E\dsv c$ and $f_E$ is as in 6.3. We thus recover a known
result from \cite{\OR, 12.17}).

\head 7. Some extensions\endhead
\subhead 7.1\endsubhead
In this subsection we assume that $W$ is irreducible and that $w\m w^\di$ is a nontrivial automorphism of $W$ 
(as a Coxeter group) such that $(w^\di)^\di=w$ for all $w\in W$. Let $\II_\di=\{w\in W;w^\di=w\i\}$ be the set of 
``$\di$-twisted" involutions of $W$. Let $M_\di$ be the free $\ca$-module with basis $(a_w)_{w\in\II_\di}$. 
Replacing $M$ by $M_\di$, $\II$ by $\II_\di$, $ws,sws$ by $ws^\di,sws^\di$ in Theorem 0.2 we obtain a true
statement. The proof is along the same lines as that in \S1,\S2: we replace $\ph'$ in $\S1$ or $\ph$ in \S2 by a 
not necessarily split Frobenius map $G@>>>G$ which induces $w\m w^\di$ on $W$ (if $G$ is not of type
$B_2,G_2$ or $F_4$) or by a Chevalley exceptional isogeny, see \cite{\CHE, 21.4,
23.7}, (if $G$ is of type $B_2,G_2$ or $F_4$ and $p=2,3,2$
respectively). Note that in this last situation, only cases (iii),(iv) appear in 0.2(a). (Indeed, in this
situation, for $s\in S$, $s,s^\di$ are not conjugate in $W$ hence $sw\ne ws^\di$ for any $w\in W$; this guarantees that the formulas in 0.2(a) specialized
with $u$ equal to an odd power of $\sqrt{p}$ involve only integer coefficients.)
The statement
obtained from Theorem 0.3 by replacing $\II$ by $\II_\di$, $\uM$ by $\uM_\di=\uca\ot_\ca M_\di$ and $P^\s_{y,w}$ 
by $P^{\s,\di}_{y,w}$ remains true with essentially the same proof.
The description of the $W$-module $M_{\di,1}=\QQ\ot_\ca M_\di$ (with $\QQ$ viewed as an $\ca$-algebra under 
$u\m 1$) given in 0.4 with $\II$ replaced by $\II_\di$ and $ws,sws$ replaced
by $ws^\di,sws^\di$) remains valid. The statement of 4.4 remains valid if $\II$ is replaced by $\II_\di$ and 
$ws,sws$ are replaced by $ws^\di,sws^\di$. (The quantities $\cm^s_{z,w}$ are defined as in 4.3 with the same 
replacements.) The analogue of 5.2(a) continues to hold. The analogues of 6.2, 6.3 continue to hold. 

\subhead 7.2\endsubhead
Theorems 0.2 and 0.3 remain valid if $W$ is replaced by an affine Weyl group. The proofs are essentially the same 
but use instead of $\cb$ the affine flag manifold associated to $G(\kk((\e)))$ (here $\e$ is an indeterminate).
We note that although $\DD$ and multiplication by $u^\nu$ (in 2.9) 
do not make sense separately, their composition $m\m u^\nu\DD(m)$ does. The results in \S4 also remain valid. 
The results in \S6 remain valid as far as $gr(M_1)$ is concerned but one cannot assert that $M_1\cong gr(M_1)$.
The generalizations of 0.2,0.3,\S4 involving $w\m w^\di$ (as in 7.1) also remain valid if $W$ is replaced by an 
affine Weyl group and $w\m w^\di$ is any automorphism of order $\le 2$ of $W$ (as a Coxeter group).

\subhead 7.3\endsubhead
Theorem 0.2(a) (and its variant for twisted involutions) remains valid if $W$ is replaced by any Coxeter group.
The main part of the proof involves the case where $W$ is a dihedral group. It is likely that Theorem
0.2(b) (and then automatically Theorem 0.3) also extend to the case of Coxeter groups.

\enddocument